\documentclass[12pt]{article} 

\def\cz{{\rm \,I \! \! \! C}}
\def\gz{{\rm Z \!\! Z}}

\def\pr{{\rm I \!  P}}
\def\r{\rightarrow}
\def\s{\subset}
\def\nz{{\rm I \!  N}}
\def\bk{{\bar\kappa}}

\def\qed{\ifmmode\sq\else{\unskip\nobreak\hfil
\penalty50\hskip1em\null\nobreak\hfil\sq
\parfillskip=0pt\finalhyphendemerits=0\endgraf}\fi}
\def\sq{\hbox{\rlap{$\sqcap$}$\sqcup$}}

\newtheorem{defi}{Definition}
\newtheorem{theo}[defi]{Theorem}
\newtheorem{prop}[defi]{Proposition}
\newtheorem{cor}[defi]{Corollary}
\newtheorem{rem}[defi]{Remark}
\newtheorem{lemma}[defi]{Lemma}

\begin{document}

\setlength{\textwidth}{6in}
\setlength{\textheight}{8.6in}
\setlength{\topmargin}{-0.3in}
\setlength{\evensidemargin}{0in}
\setlength{\oddsidemargin}{0in}

\title{Logarithmic Surfaces and Hyperbolicity}
\author{Gerd Dethloff and
Steven S.-Y. Lu
\thanks{Partially supported by an NSERC grant and the DFG-Leibniz 
program.}}
\date{}
\maketitle

\section{Introduction and results}
We first introduce some notations (see
section 2 for the precise definitions).
Let $\bar X$ be a projective manifold
and $D \s \bar X$  a normal crossing divisor. We call
the pair
$(\bar X , D)$ a log manifold and denote $X = \bar X \setminus D$.
It is called a log surface if $\dim \bar X=2$.
Let $T^*_{\bar X}$ be its cotangent bundle and
$\bar T ^*_{X}$ its log cotangent bundle.  We denote by
$q_{\bar X} = \dim  H^0(\bar X, T^*_{\bar X})$ its irregularity
and
$\bar q_{ X} = \dim  H^0(\bar X, \bar T^*_{X})$ its log
irregularity.
We denote its log canonical bundle  by
$\bar K_{X} = \bigwedge^{\dim \bar X} \bar T^*_{X} = K_{\bar
X}( D)$
and its log Kodaira dimension
by  $\bar \kappa_{X} = \kappa (\bar X, \bar K_{X})$, the $L$-dimension
of $\bar K_{X}$. We call $(\bar X,D)$ to be of log general type if
$\bar \kappa_{X}= \dim \bar X$.
Finally let $\alpha_X : X \r {\cal A}_X$ be the quasi-Albanese map.
It is a holomorphic map which extends to a rational map
$\bar \alpha_X : \bar X - \r \bar {\cal A}_X$ (Iitaka '76 \cite{Ii2}),
where $\bar {\cal A}_X$ is some compactification of ${\cal A}_X$ (see  \S 2.2).

We know by the log-Bloch theorem (Noguchi '81
\cite{No2}) that for any log manifold such that $\bar q_X >{\rm
dim}\bar X $,
any entire holomorphic curve $f:\cz \r X$ is algebraically degenerate
(this means $f(\cz)$ is contained
in a proper algebraic subvariety of $\bar X$). More generally, by
results of Noguchi-Winkelmann '02 \cite{NW} one has (defining the
log structure
and $\bar q_{X}$ on a K\"{a}h\-ler manifold in a similar way as above):
\begin{theo} [Noguchi, Noguchi-Winkelmann] \label{t1}Let $\bar X$ be a
compact
K\"{a}h\-ler manifold and $D$ be a hypersurface in $\bar X$. If $\bar
q_{X} > \dim  \bar X$,
then  any entire holomorphic curve $f:\cz \r X$ is analytically
degenerate (this means contained
in a proper analytic subset of $\bar X$).
\end{theo}
In this paper we are interested
in the case of log surfaces $(\bar X , D)$ with $\bar q_X =\dim \bar X = 2$.\\

The first part of this paper deals with the case of surfaces of log general type, that is 
$\bar \kappa_{X}= 2$. We restrict our attention to Brody curves, this means 
 entire curves
$f:\cz \r X$ with bounded derivative
$f'$  in $\bar X$, for which we resolve the problem
completely in the following Main Theorem of our paper.

\begin{theo} \label{t2} Let $(\bar X,D)$ be a log surface  with log
irregularity $\bar q_X=2$ and with log Kodaira dimension $\bar
\kappa_{X}=2$.
    Then every Brody curve $f:\cz \r X$ is algebraically
degenerate.
\end{theo}

\noindent The proof and some applications are given in section 3. \\

The second part of this paper complements Theorem~\ref{t2}.

Let $(\bar X,D)$ be a log surface  with log irregularity $\bar q_X=2.$
Let $\bar \alpha_X : \bar X - \r \bar {\cal A}_X$ be the compactified quasi-Albanese map, $I$ its
finite set of points of indeterminacy and 
$\bar \alpha_0=\bar \alpha_X|_{\bar X\setminus I}$.
In the case of dominant $\bar \alpha_X$, we consider the following condition:
\begin{itemize}
\item[\rm{(*)}] For all $z\in {\cal A}_X$ and $\bar E$ a connected component of
the Zariski closure of $\bar \alpha_0^{-1}(z)$ with $\bar E\cap X\neq \emptyset$,
any connected component of $D$ intersecting $\bar E$ is contained in $\bar E$ 
(i.e. $\bar E$ is a connected component of $\bar E\cup D$).
\end{itemize}

\noindent We remark that condition (*) can be expressed intrinsically (see subsection 4.2)
and is implied by the condition that all the fibers of $\alpha_X:X\r {\cal A}_X$ are compact.
In particular, this condition is much weaker than the properness of $\alpha_X$.

In the case $\bar \kappa_{X}= 1$ we have the following result:

\begin{theo} \label{neut} Let $(\bar X,D)$ be a log surface  
with log irregularity $\bar q_X=2$ 
and with log Kodaira dimension $\bar \kappa_{X} =1$. 
Assume condition {\rm (*)} in the case of dominant $\bar \alpha_X$.
Then every entire curve $f:\cz \r X$ is algebraically
degenerate.
\end{theo}

\noindent
 The proof and some applications are given in section 4. \\

As our counterexample
in Proposition~\ref{p6} shows, the additional condition (*) is necessary for 
the theorem to hold. Some reflections on our proof will reveal also that our 
condition (*) can in fact not be weakened further, at least when mild restrictions are imposed,
for whose discussion and generalization to higher dimensions will be 
relegated to another paper.

\begin{rem}\label{r7}
It is easily obtained from Hodge theory due to Deligne '71 \cite{D}
(see, for example, Catanese '84 \cite{C})
that we have, for $(\bar X, D)$ a log surface:
\begin{equation}\label{1}{\rm rank}_{\gz}{\rm NS}(\bar X) \geq
{\rm rank}_{\gz}\{c_1(D_i)\}_{i=1}^k = k - \bar q_{X} + q_{\bar X}
\end{equation}
where $D_1,...,D_k$
are the irreducible components of $D$ and NS$(\bar X)$ denotes the
Neron-Severi group
of $\bar X$. This
may be deduced from  the proof of Theorem
1.2(i) of Noguchi-Winkelmann
'02 \cite{NW}, p.605.  But there does not seem to be an easy way to
profit from this, unless one assumes some bound on the Neron-Severi
group of $\bar X$.
\end{rem}

\vspace{-3mm}
Connected with this, we would like to mention again the work of
Noguchi-Winkelmann '02
\cite{NW}, which deals with the question of algebraic degeneracy 
in all dimensions and even with K\"{a}hler  manifolds,
especially with log tori or with log manifolds having small Neron-Severi
groups, under the additional condition
that all irreducible components $D_i$ of $D$ are ample. But as can be
seen from the
equation (\ref{1}) above,  their results never concern
    the case of log surfaces with log irregularity $\bar q_{X} \leq 2$.
    
We would also like to mention the preprint of Noguchi-Winkelmann-Yamnoi '05
\cite{NWY2} that we just received, which deals with the case of algebraic manifolds
whose quasi-Albanese map is proper. 
More precisely, Noguchi-Winkelmann-Yamanoi '05
\cite{NWY2}  deals with arbitrary holomorphic curves into arbitrary dimensional
algebraic varieties of general type (this is the essential case), but only with proper Albanese map.
Although we deal entirely with surfaces and mainly with Brody curves,
our result do not require any condition
in the case of log general type and our condition (*), which we use in the
general case, is still much weaker than the condition of properness of the
quasi-Albanese map as can, for example,  be seen from the following simple example:
Let $X$ be the complement of a smooth ample divisor $D$ in a abelian
surface $\bar X$. Then $\bar K_X=D$ and so $\bar \kappa_X=2$. Also
$\bar q_X=q_{\bar X}=2$ by equation (\ref{1}). This means that the compactified
quasi-Albanese map is the identity map and so the Albanese map is not
proper.\\

We now give an indication of our methods of proof.

We first discuss the ideas of the proof for Theorem~\ref{t2}. 
We will  reduce the proof of this theorem by a result of McQuillen and ElGoul 
and by log-Bloch's theorem to the 
claim  that under the conditions of Theorem~\ref{t2}, $\alpha_X\circ f: \cz \r {\cal A}_X$ is a translate of a complex
one parameter subgroup of ${\cal A}_X$.\\
In the case $q_{\bar X}=2$, the compactified quasi-Albanese map $\bar
\alpha_{X}$
is a morphism and so the claim is trivial.
If $q_{\bar X}<2$, $\bar \alpha_{X}$ can have points of indeterminacy
so that Brody curves are not preserved by $\bar \alpha_{X}$ in general.
But from value distribution theory, the order of growth of a
holomorphic curve is
preserved under rational maps and Brody curves are of order at most two.
Using this, the key analysis in this proof consists of a detailed study
of the geometry
of the quasi-Albanese map (in particular at its points of
indeterminacy) with respect to $f$
to reduce the order of  $\bar \alpha_{X} \circ f$ to one or less.
Then $\bar \alpha_{X} \circ f$ is either constant or a leaf
of a linear foliation on ${\cal A}_{X}$. 
We do this componentwise where in the case $q_{\bar X}=1$, we use the
fact from
Noguchi-Winkelmann-Yamanoi '02 \cite{NWY} that one can choose a metric
on $\bar {\cal A}_{X}$
which lifts to the product metric on the universal cover $\cz \times
\pr_{1}$ of
$\bar {\cal A}_{X}$.
In the case $q_{\bar X}=0$ we take rational monomials of the components
of $\bar \alpha_{X}$ motivated by
arranging residues in a way that allows us to control the points of
indeterminacy of the
resulting map with respect to $f$.

We now discuss our proof for Theorem~\ref{neut}.
We first prove the analogue in Proposition~\ref{neutt} 
of the structure theorem of Kawamata for open
subsets of finite branched covers of semi-abelian varieties. We follow essentially
the original ideas of Kawamata but with several new ingredients. For example,
in the case $\bar \kappa_X=1$, one needs to observe that 
even though the quasi-Albanese map to the semi-abelian 
variety is not proper, the restriction to the generic fiber is. In the case $\bar\kappa_X=0$,
we need to observe
that a complement of a (singular) curve in a semi-abelian variety
is of log general type unless the curve is a translate of a algebraic subgroup. 
We reduce this observation by the addition theorem of Kawamata to the case
where the semi-abelian variety is a simple abelian variety. For
Theorem~\ref{neut}, the main observation is that condition (*) 
is equivalent to a condition on the Stein factorization of a desingularization of
the quasi-Albanese map and that this allows us to use
Proposition~\ref{neutt} to conclude that the base of the Iitaka fibration is hyperbolic.\\

We  remark that we can give an elementary proof of the
result of McQuillan and
El Goul in the case of linear foliations on ${\cal A}_{X}$ by using
techniques similar
   to those given in Bertheloot-Duval '01 \cite{BD}.  \\
   
The first named author would like to thank C.Peters and F.Catanese
for valuable discussions some years ago on configurations
of boundary curves of log surfaces in relation to log 1-forms.
Finally we would like to  thank the referee for many suggestions
on how to clarify and improve the presentation of this paper.

\newpage

\section{Some Preliminaries }

\subsection{Log manifolds and residues of log 1-forms}

Let $\bar{X}$ be a complex manifold with a normal
crossing divisor $D$. This means that around any point $x$ of $\bar{X}$,
there exist local coordinates $z_1,...,z_n$
centered at $x$ such that $D$ is defined by
$z_1z_2...z_l=0$ in some neighborhood of $x$ and for some $0 \leq l\leq
n$.
The pair $(\bar{X}, D)$ will be called a {\it log-manifold}.
Let $X= \bar{X} \setminus D$.

Following Iitaka '82 \cite{Ii}, we define the logarithmic
cotangent sheaf
   $$\bar{T}^{*}_{X} = \bar{\Omega}X
= \Omega(\bar X, {\rm log} D)$$  as the locally free
subsheaf of the sheaf
of meromorphic 1-forms on $\bar{X}$,
    whose restriction to $X$ is $T^{*}_{X}= \Omega X$
(where we identify from now on vector bundles and their sheaves of sections)
and whose localization at $x\in \bar X$ is of the form
\begin{equation} \label{9}
(\bar{T}_{X}^{*})_x=
\sum_{i=1}^l {\cal O}_{\bar X,x}{dz_i\over z_i} +
\sum_{j=l+1}^n {\cal O}_{\bar X,x} dz_j,
\end{equation}
where the local coordinates $z_1,...,z_n$ around $x$ are chosen as
before.
Its dual, the logarithmic tangent sheaf $\bar{T}_{X} =T (\bar X, - {\rm
log} D)$,
is a locally free subsheaf of
the holomorphic tangent bundle $T_{\bar{X}}$ over $\bar{X}$.
Its restriction to $X$ is identical to $T_{X}$,
and its localization at  $x\in \bar X$ is of the form
\begin{equation} \label{10}
(\bar{T}_{X})_x=
\sum_{i=1}^l {\cal O}_{\bar X,x}z_i{\partial \over \partial z_i} +
\sum_{j=l+1}^n {\cal O}_{\bar X,x} {\partial \over \partial z_j}.
\end{equation}

\noindent 
Let $\omega $ be a log 1-form defined around $x$, so that by (\ref{9})
we have
\begin{equation} \label{r}
{\omega}_{x}=
\sum_{i=1}^l (h_{i})_{x}{dz_i\over z_i} +
\sum_{j=l+1}^n (h_{j})_{x} dz_j
\end{equation}
Then we call, for $i=1,...,l$, the complex number $(h_{i})_{x}(x) \in
\cz$
the {\it residue} of $\omega$ at $x$ on the local irreducible branch
of $D$
given by $z_{i}=0$.
Since we do not assume simple normal crossing,
we may have several such local irreducible components for $D$ at $x$
even if $D$ is irreducible. But if $\bar X$ is compact, it is easy to
see that for any
(global) irreducible component $D_{j}$ of $D$, the residue is constant
on $D_{j}$
(meaning it is the same for all points $x \in D_{j}$). In fact, whether
$\bar X$
is compact or not, we have the exact sequence of sheaves:
\begin{equation}\label{2!}
0 \r \Omega_{\bar X} \r \bar \Omega_{X} \stackrel{\rm Res}{\r}
{\cal O}_{\hat D} \r 0\:,
\end{equation}
where $\hat D$ is the normalization of $D$.\\

\subsection{Quasi-Albanese maps}

We first recall the definition and some basic facts on
semi-abelian varieties (see  Iitaka '76 \cite{Ii2}).

A quasi-projective variety $G$ is called a {\it semi-abelian variety}
if it is a complex commutative Lie group which admits an
exact sequence of groups
\begin{equation}\label{2}
0 \r (\cz^*)^{l} \r G \stackrel{\pi}{\r} A \r 0\:,
\end{equation}
where $A$ is an abelian variety of dimension $ m$. An important
point in our analysis is that this exact sequence is not unique,
but depends on a choice of $l$ generators for the kernel of $\pi$, which is an
algebraic torus of dimension $l$.

   From the standard compactification  $(\cz^*)^{l} \s (\pr_1)^l$, which
is
   equivariant with respect to the  $(\cz^*)^l$ action,
   we obtain a  completion
$\bar{G} $ of $G$ as the $(\pr_1)^l$ fiber bundle associated to
the $(\cz^*)^l$ principal bundle $G\r A$.  This is
a smooth compactification of $G$  with a simple normal
crossing boundary divisor
$S$.
The projection map
$$\bar \pi : \bar G \r A$$
    has the structure of a $(\pr_{1})^{l}$-bundle. Here, $\bar \pi$ and $ \bar G$
  depend on the choice of the $l$ generators that identified our algebraic torus
  as $(\cz^*)^{l} $.

We denote the natural action of $G$ on
$\bar{G}$ on the right as addition.
It follows that the exponential map from the Lie algebra $\cz^n$ to $G$
is a group homomorphism and, hence, it is also the universal covering
map of
$G=\cz^n / \Lambda$, where
    $\Lambda=\pi_1(G)$
is a discrete subgroup of $\cz^n$ and $n= m+l$.

Following Iitaka '76 \cite{Ii2}, we have the following
explicit trivialization of
the log tangent and cotangent bundles of $\bar{G}$:
Let $z_1, ...,z_n$ be the
standard coordinates of $\cz^n$. Since $dz_1,...,dz_n$ are invariant
under
the group action of translation on $\cz^n$, they descend to
forms on $G$. There they extend to
logarithmic forms on $\bar G$ along $S$, which are elements of
$H^0(\bar G, \bar T_{G}^{*})$. These logarithmic $1$-forms
are everywhere
linearly independent on $\bar G$.
Thus, they
globally trivialize the vector bundle
$\bar T_{G}$.
Finally, we note that these log 1-forms
are invariant under the group action of $G$ on $\bar{G}$, and, hence,
    the associated trivialization of $\bar T_{G}$ over $\bar G$
is also invariant.\\

Let now $(\bar X, D)$ be again a log surface and
    $\alpha_{\bar X}: \bar X \r {\cal A}_{\bar X}$ the Albanese map of
$\bar X$
    (it can be constant if $q_{\bar X}=0$). Taking into account also the
log 1-forms,
    Iitaka '76 \cite{Ii2} introduced the quasi-Albanese map $\alpha_{X}:X
\r {\cal A}_{X}$,
    which is a holomorphic map to the semi-abelian variety ${\cal 
A}_{X}$,
which comes equipped with the exact sequence
    \begin{equation}\label{3}
0 \r (\cz^*)^{l} \r {\cal A}_{X} \stackrel{\pi}{\r} {\cal A}_{\bar X}
\r 0
\end{equation}
(Iitaka makes a noncanonical choice of $l$ generators for the algebraic torus for
this construction).
We have the commutative diagram
\begin{equation}\label{4}
\begin{array}{ccc}
X
& \stackrel{ \alpha_X}{ \r }
& {\cal A}_X \\
    & & \\
& \stackrel{ \alpha_{\bar X}}{ \searrow }
& \downarrow \pi\\
    & & \\
&
& {\cal A}_{\bar X}
    \end{array}
\end{equation}
Iitaka also proved that $\alpha_{X}$ extends to a rational map
$\bar \alpha_{X}: \bar X - \r \bar {\cal A}_{X}$, and the diagram
(\ref{4}) extends to
\begin{equation}\label{5}
\begin{array}{ccc}
\bar X
& - \stackrel{ \bar \alpha_X}{ \r }
& \bar {\cal A}_X \\
    & & \\
& \stackrel{ \alpha_{\bar X}}{ \searrow }
& \downarrow  \bar \pi\\
    & & \\
&
& {\cal A}_{\bar X}
    \end{array}
\end{equation}

    In general the $(\pr_{1})^{l}$-bundle
$\bar {\cal A}_{X} \stackrel{\bar \pi}{\r} {\cal A}_{\bar X}$ is not
trivial.
But Noguchi-Winkelmann-Yamanoi '02 \cite{NWY} observed
that the transition functions of the $(\pr_1)^{l}$-bundle
$\bar \pi :\bar {\cal A}_X \r {\cal A}_{\bar X}$
(as the structure group $ (\cz^*)^{l}$ can always be reduced to
the subgroup defined by $|z_i|=1,\ i=1,..,l$)
can be chosen to be
isometries with respect to
    the product Fubini-Study metric on $(\pr_1)^{l}$.
\begin{prop} \label{p9}There exists a metric $h$ on $\bar {\cal A}_X$
so that the universal
cover map $$\left(\cz^{m} \times (\pr_1)^{l}, {\rm eucl.} \times {\rm
product\ FS}\right) \r (\bar {\cal A}_X, h)$$ is a local isometry.
\end{prop}

 We first consider the case $\dim \bar X
=2$ and $q_{\bar X}=0$. This means we have a morphism $\alpha_X: X \rightarrow 
\cz^* \times \cz^*$ extending to a rational map 
$\bar \alpha_X: \bar X  - \rightarrow 
\pr_1 \times \pr_1$.
Let
$\omega \in H^{0}(\bar X, \bar T^{*}_{X})$ be a log 1-form with residues
$a_{j} \in \gz$ along the irreducible components $D_{j}$ of $D$,
$j=1,...,k$.
We define a  holomorphic function
\begin{equation}\label{6a}
\Phi : X \r \cz^{*}; \; \Phi (x)= \exp \left(\int_{x_{0}}^{x}\omega
\right),
\end{equation}
where $x_{0}$ is a fixed point in $X$. It is well defined since by
the condition $q_{\bar X}=0$, there are no non trivial periods in $\bar X$,
and since the periods around the components of the divisor $D \subset \bar X$
of the integral have values which are entire multiples of $2 \pi$, and, hence,
are eliminated by taking the $\exp$-function. 
We claim that this function extends to a rational function $\bar \Phi :
\bar X -\r \pr_{1}$:
Let  $P \in D_1$ (respectively $P \in D_1 \cap D_2$) and let
$z_1,z_2$ be local coordinates
around $P$ such that
   $D_1= \{z_1=0 \}$ (respectively $D_1= \{z_1=0 \}$
and $D_2= \{z_2=0 \}$).
Then it follows from equation (\ref{6a}) that there exists a
holomorphic function
$h:U(P) \r \cz^*$ on a neighborhood $U(P)$ of $P$ such that
\begin{equation}\label{8}
    \bar \Phi(z_1,z_2) = z_1^{a_{1}}h(z_1,z_2)\;\;\;\left( {\rm
respectively}\;
\bar \Phi(z_1,z_2) = z_1^{a_{1}}z_2^{a_{2}}h(z_1,z_2)\right).
\end{equation}
In more detail, if $D_1= \{z_1=0 \}$ around $P$, then we have
$$\omega = h_{1}(z){dz_1\over z_1} + h_{2}(z)dz_{2} = a_{1}{dz_1\over
z_1} +
\frac{h_{1}(z)-a_{1}}{z_{1}}dz_{1}+ h_{2}(z)dz_{2}\, . $$
By the first Riemann extension theorem, the function
$\frac{h_{1}(z)-a_{1}}{z_{1}}$
extends to a holomorphic function, so we have
$\omega =a_{1}{dz_1\over z_1} + \omega_{hol} $, where $\omega_{hol}$ is a
holomorphic form.
Now by (\ref{6a}) we get
$$\Phi (z)= \exp \left(\int_{x_{0}}^{z}\omega \right)=
\exp \left(\int_{x_{0}}^{z}a_{1}{dz_1\over z_1} + \omega_{hol}\right) =
\exp \left( a_{1} (\log z_{1}+2\pi i \gz)\right) \cdot \exp (\tilde h(z))$$
$$
=\left( \exp  \left( \log z_{1}+2\pi i \gz\right)\right)^{a_{1}}
   \cdot \exp \left( \tilde h(z)\right)
=z_{1}^{a_{1}}h(z)$$
The other equality follows in the same way.

   From this we get by the second Riemann extension theorem, GAGA and by
the local description
of points of indeterminacy the following.

\begin{prop} \label{p12}
The  holomorphic map $\Phi : X \r \cz^{*}$ given by  (\ref{6a}) extends to a rational
function
     $\bar \Phi: \bar X \r \pr_1$. It is a morphism outside the points of
intersection of pairs of different
irreducible components $D_1, ..., D_k$ of $D$. In particular points of indeterminacy
never occur at
self-intersection points of a component.
More precisely, a
point of $D_{j_{1}}\cap D_{j_{2}} \subset \bar X$ is in the set $I$ of  points
of indeterminacy of $\bar \Phi$  iff
$a_{j_{1}} \cdot a_{j_{2}} < 0$, and $\bar \Phi(D_j\setminus I)\equiv 0$ (respectively 
$\infty$) iff the residue $a_j$ of $\omega$ along $D_j$ is $> 0$ (respectively
$<0$).
\end{prop}

Using this we now describe the
components
of the map $\alpha_{X}=\left((\alpha_{X})_{1}, (\alpha_{X})_{2}\right):
X \r
(\cz^{*})^{2}$
in more detail.  The following basic facts
follow from the exact sequence $(\ref{2!})$  and duality in
Hodge theory,
and they can be found
for example in Noguchi-Winkelmann '02 \cite {NW}:

\begin{prop} \label{p11} We can choose the basis
$\omega_1,\omega_2 \in
H^0(\bar X, \bar T^*_{X})$ such that the residue $a_{ij} \in \gz$
for  all $ i=1,2$ and $j=1,...,k$,
where $a_{ij}$ is the residue of $\omega_i$
along the irreducible component $D_j$ of $D$.
    The matrix of residues so obtained  has rank 2:
\begin{equation}\label{7}
\begin{array}{ccccc}
&
D_1&
D_2 &
... &
D_k \\
\omega_1 &
a_{11} &
a_{12}  &
... &
a_{1k} \\
\omega_2 &
a_{21} &
a_{22}  &
... &
a_{2k}
\end{array}
\end{equation}
\noindent
and we have, as in \rm{(\ref{6a})},
\begin{equation}\label{6}
(\alpha_{X})_{i}(x)= \exp \left(\int_{x_{0}}^{x}\omega_{i}\right),
\end{equation}
where $x_{0}$ is a fixed point in $X$.
\end{prop}

Henceforth, in the case $q_{\bar X}=0 $, we assume that the components of $\alpha_X$ are
given by such a choice of basis.\\

Next we consider the case where $\dim \bar X=2$ and $ q_{\bar X}=1$. 
Let $x\in D$ be a point. By diagram (\ref{5}) and the notation thereof, as 
${\cal A}_{\bar X}$ is an elliptic curve, there is a small open neighborhood $W$
of $\alpha_{\bar X}(x)$ such that $\bar\pi^{-1}(W)\simeq \pr_1\times W$.
Let $\bar V\s \alpha_{\bar X}^{-1}(W)$ be a small open ball 
centered at $x$ in $ \bar X$. Then on $V=\bar V\cap X$,
$\alpha_X$ can be written as 
$$\alpha_X=(\alpha_{\bar X},\Phi)$$
where $\Phi : V\r \cz^*$ is as in (\ref{6a}) (see Noguchi-Winkelmann '02 \cite{NW}).
The same argument as we gave for Proposition~\ref{p12} gives us:

\begin{prop} \label{p12a}
The  holomorphic map $\Phi : V \r \cz^{*}$ given by  (\ref{6a}) extends to a rational
function
     $\bar \Phi: \bar V \r \pr_1$. It is a morphism outside the points of
intersection of pairs of different
irreducible components $D_1, ..., D_k$ of $D$. In particular points of indeterminacy
never occur at
self-intersection points of a component.
More precisely, a
point of $D_{j_{1}}\cap D_{j_{2}} \subset \bar X$ is in the set $I$ of  points
of indeterminacy of $\bar \Phi$  iff
$a_{j_{1}} \cdot a_{j_{2}} < 0$, and $\bar \Phi(D_j\setminus I)\equiv 0$ (respectively 
$\infty$) iff the residue $a_j$ of $\omega$ along $D_j$ is $> 0$ (respectively
$<0$).
\end{prop}

\subsection{Brody curves, maps of order 2 and limit sets of entire
curves}

Let $(\bar X, D)$ be a log manifold and $f:\cz \r X$ be an entire
curve. We recall that
$f$ is a Brody curve if the derivative of $f$ with respect to some (and
so any) hermitian
metric on $\bar X$ is bounded.

Following Noguchi-Ochiai '90 \cite{NO}, we have the characteristic
function
$T_{f}(r,\omega)= \int_{1}^{r}\frac{dt}{t}\int_{|z|<t}f^{*}\omega$ of a
holomorphic map
$f:\cz \r Y$ with respect to a real continuous $(1,1)$-form $\omega$
on a K\"{a}hler manifold $Y$. If $Y$ is compact and $\omega_{H}$
denotes the $(1,1)$-form
associated to a hermitian metric $H$ on $Y$, then it is easy to see
(\cite{NO}, (5.2.19)) that
$\rho_{f}:=\overline{\lim}_{r \r \infty}\frac{\log T_{f}(r,
\omega_{H})}{\log r}$
is independent of the hermitian metric $H$.
The map $f$ is
said to be of {\it order at most 2} if
    $\rho_{f} \leq 2$.
    Since the derivative of a Brody curve in the projective variety $\bar
X$ is  bounded
    with respect to  hermitian metrics on $\bar X$, Brody curves
    are easily seen to be of order at most 2. By a classical 
theorem of Weierstrass, we get also
    that a curve $f:\cz \r \pr_{N}$ of order at most 2 which omits the
coordinate hyperplanes
    can be written in the form $\big(1:\exp (P_1(z)):...:\exp (P_n(z))\big)$ where
the $P_{i}$'s are polynomials  of degree at most 2 in the variable
    $z \in \cz$.

    We now prove that the property of having  order at most 2 is 
preserved
under rational maps
    (see also \cite{DSW1} for a similar result).
    \begin{lemma}\label{l13}
    Let $f:\cz \r \pr_{N}$ be a curve of order at most 2 and $$R:\pr_{N}
-\r \pr_{M}$$ be
    a rational map (not necessarily dominant) such that $f(\cz)$ is not
contained in the set of indeterminacy of $R$.
     Then the curve $R \circ f : \cz \r \pr_{M}$ is of order at most 2.
    \end{lemma}
    \noindent {\bf Proof} Let $f=(f_{0}:...:f_{N})$ be a reduced
representation and
    $$R=(R_{0}:...:R_{M})$$ be a (not necessarily reduced) representation
by polynomials $R_{0},..., R_{M}$ of degree $p$.
    Then $$\left(f_{0}^{p}:...:f_{N}^{p}:R_{0}\circ f:...:R_{M}\circ
f\right)$$ is a
reduced representation
    of a curve  $F: \cz \r \pr_{N+M+1}$, and without loss of generality
$R_{0} \circ f \not\equiv 0$.
    We have by \cite{NO}, p.183
    $$T_{F}(r, \omega_{{\rm FS}})=\int_{0}^{2 \pi}\frac{1}{2}\log
\left( |f_{0}^{p}|^{2}+...+|f_{N}^{p}|^{2}+
    |R_{0} \circ f|^{2}+...+ |R_{M} \circ f|^{2}\right) d\theta +O(1)\leq
$$
    $$\int_{0}^{2 \pi}\frac{1}{2}\log
\left( |C|(|f_{0}^{p}|^{2}+...+|f_{N}^{p}|^{2})\right) d \theta +O(1)
    =\int_{0}^{2 \pi}\frac{p}{2}\log \left(
|f_{0}|^{2}+...+|f_{N}|^{2}\right) d \theta
+O(1)$$
$$=p \cdot T_{f}(r, \omega_{{\rm FS}}) +O(1) $$
    By Noguchi-Ochiai '90 \cite{NO} (5.2.29) and (5.2.30) we get
    $$T_{R \circ f}(r, \omega_{{\rm FS}}) \leq
    \sum_{j=1}^{M}T_{\frac{R_{j }\circ f}{R_{0}\circ f}}(r, \omega_{{\rm
FS}})
    \leq M \cdot  T_{F}(r, \omega_{{\rm FS}}) +O(1) \leq M \cdot p \cdot
T_{f}(r, \omega_{{\rm FS}})+O(1)$$
    From this the lemma follows. \qed
    \bigskip
    Finally, we need the definition and some simple observations on the limit 
set
of an entire curve
    $f: \cz \r \bar X$  as given by Nishino-Suzuki '80 \cite{NS}
    (Proposition 1 on p.463 and Proposition 3 on p.466).  
    For $r> 0$, put
    $\Delta_{r}^c= \{z \in \cz : |z|>r \}$. Let
$\overline{f(\Delta_{r}^c)}
\subset \bar X$
    be the closure (with respect to the usual topology) of
$f(\Delta_{r}^c)$
in $\bar X$,
    and $f(\infty):=\bigcap_{r>0}\overline{f(\Delta_{r}^c)}$. We remark
that
$f(\infty)$
    is exactly the set of all points $p \in \bar X$ such that there 
exists
a sequence
    $(z_{v})_{v \in \nz}$ with $|z_{v}|\r \infty$ and $f(z_{v}) \r p$.

    \begin{prop}\label{p14} Let $f: \cz \r X$ be as above
    for a log surface $(\bar X, D)$.
    \begin{itemize}
    \item[a)] $f$ extends to a holomorphic map $\hat{f}: \pr_{1} \r \bar
X$ iff $f(\infty)$ consists
    of exactly one point.
    \item[b]) If $f(\infty)$ is contained in a proper algebraic subvariety
$C \subset \bar X$
    and is not a point,
    then $f(\infty)$ is equal to a    union of some of the irreducible
components of $C$.
    \end{itemize}
    \end{prop}

    \begin{lemma}\label{limit1}
For $a\in \cz^*$, consider the entire curve 
$f: \cz\r \cz^*\times\cz^*,$ given by $z\mapsto (\exp (az), \exp (z^2))$. 
Then $f(\infty)\subset \pr_1\times \pr_1$  contains $\{0,\infty\}\times \pr_1$.

\end{lemma}

\noindent{\bf Proof }  Since $f(\infty)$ is closed and since we can change the
sequence $(z_\mu)_{\mu \in \nz}$ to the sequence $(-z_\mu)_{\mu\in \nz}$, it suffices
to prove that $f(\infty)$ contains the points $(\infty, c)$ with $c\in \cz^*$. Let  $c \in \cz^*$ fixed and
$b\in \cz$ be such that $c=\exp b$. We consider the four sequences 
$z_\mu=\pm \sqrt{b \pm 2\pi i\mu}$, $\mu \in \nz$, where the four different sequences are
obtained by the four different choices of signs. Then we always have $\exp (z_{\mu}^2)=c$,
and for the arguments we have that $\arg (z_{\mu})$ converges to one of the numbers
$\frac{\pi}{4}, \frac{3\pi}{4},\frac{5\pi}{4},\frac{7\pi}{4}$ according to the choice of signs.
Then $|az_{\mu}| \rightarrow \infty$ and $\arg (az_{\mu})$ converges to one of the numbers
$\arg (a)+\frac{\pi}{4}, \arg (a)+\frac{3\pi}{4}, \arg (a)+\frac{5\pi}{4}, \arg (a)+\frac{7\pi}{4}$. 
Since one of these four numbers is contained, modulo $2\pi i \gz$, in the open interval $(-\frac{\pi}{2}, \frac{\pi}{2})$, there exists a choice of signs such that $\exp (az_\mu) \rightarrow \infty$. \qed

    \begin{lemma}\label{limit2}
Let $(\bar X,D)$ and $(\bar Y, E)$ be log surfaces, 
$\Psi: X \r Y$ a morphism and $\bar \Psi: \bar X - \r \bar Y$ the 
rational extension of $\Psi$. Then $$\bar \Psi (f(\infty)) \supset \Psi\circ f(\infty),$$
where we understand $\bar \Psi(x)$ to be a curve for $x$ a point of inderminancy of $\bar \Psi$.
\end{lemma}

\noindent{\bf Proof } Let $Q\in \Psi\circ f(\infty)$. Then there is a sequence 
$(z_\mu)_{\mu \in \nz}$ in $\cz$ such that $\lim_{\mu\r \infty} |z_\mu| =\infty$
and $\lim_{\mu\r \infty} \Psi\circ f(z_\mu) =Q$. Let $\bar{\hat\Psi}: \bar{\hat X} \r \bar Y$
be a desingularization of $\bar \Psi$ via $\bar r:   \bar{\hat X} \r \bar X$   and 
$\hat f: \cz\r \bar{\hat X}$ the lift of $f$ to $\bar{\hat X}$. 
As $\bar{\hat X}$ is compact, after passing to a subsequence,
$\hat f (z_\mu)$ converges to a point $\hat P\in \bar{\hat X}$ and by continuity
of $\bar r$, $P=\bar r(\hat P)$ lies in $f(\infty)$. Hence, the continuity 
of $\bar{\hat\Psi}$ gives $Q=\bar{\hat\Psi}(\hat P)\in 
\bar \Psi(P)\subset \bar\Psi(f(\infty))$. 
\qed

\subsection{Foliations and the theorem of McQuillan-ElGoul}
We first recall some basic notations and facts on foliations as given in
Brunella '00 \cite{Br}: A foliation on the surface $\bar X$ can be 
defined by a collection of 1-forms 
$\omega_i \in \Omega^1_{\bar X}(U_i)$ with isolated singularities 
such that
$U_i,i\in I$ is an open covering of $\bar X$, and that we have 
$\omega_i = f_{ij} \cdot \omega_j$ on $U_i \cap U_j$ with 
$f_{ij} \in {\cal O}_{\bar X}^*(U_i \cap U_j)$. These isolated singularities
are called singularities of the foliation. The local integral curves of the
forms $\omega_j$ glue together, up to reparametization, giving the so-called
leaves of the foliation.
Any meromorphic 1-form $\omega$ on $\bar X$ gives a foliation, namely 
choose an open covering  $U_i,i\in I$ of $\bar X$ and meromorphic
functions $f_i$ on $U_i$ such that $\omega_i:=f_i \cdot \omega|_{U_{i}}$ 
are holomorphic
forms on $U_i$ with isolated singularities, only. Then $\omega_i, i \in I$
gives a foliation. 
By Brunella '00 \cite{Br} for algebraic $\bar X$ any foliation can be
obtained like this, and, moreover, there is a one to one correspondance
between foliations on the one hand and rational 1-forms modulo 
rational functions on the other hand.

 We  will use the following result of
McQuillan '98 \cite{MQ} and Brunella '99 \cite{Br0}, extended to the log
context
by El Goul '03 \cite[Theorem 2.4.2]{EG} (in order to avoid misunderstandings
we would like to point out that in this Theorem 2.4.2 of El Goul '03 
\cite{EG} the foliation does 
{not} need to be tangent to the boundary divisor):
\begin{theo}\label{t15} Let $(\bar X, D)$ be a log surface of
log
general type. Let $f:\cz \r X$ be an entire curve. Suppose that there
exists a foliation
${\cal F}$ on $\bar X$ such that $f$ is (contained in) a leaf of ${\cal
F}$.
Then $f$ is algebraically degenerate in $\bar X$.
\end{theo}

What we will need for  the proof of Theorem \ref{t2} is only a corollary
of a very special
case of Theorem \ref{t15}:
\begin{prop}\label{p16}
Let $(\bar X,D)$ be a log surface of log general type. Let\\
   $\Psi : X \r {\cal A}_{X}$ be a dominant
morphism which
extends to a rational map  $\bar \Psi : \bar X - \r \bar {\cal A}_{X}$.
Let $f: \cz \r X$
be an entire curve. Assume that the map
$\Psi \circ f:\cz \r {\cal A}_{X}$ is linearly degenerate
with respect to the universal cover $\cz^{2} \r {\cal A}_{X}$.
    Then $f$ is
algebraically degenerate.
\end{prop}

\noindent {\bf Proof }
The map $ \Psi \circ f : \cz \r {\cal A}_{X}$ is linearly
degenerate with respect to the universal cover $\cz^{2} \r {\cal
A}_{X}$.
Hence, $\Psi \circ f$ is a leaf of the linear foliation given by a
1-form with
constant coefficients on $\cz^{2}$, descending to a nowhere vanishing
log 1-form
$\omega$
on $\bar {\cal A}_{X}$ (see subsection 2.2). It corresponds to a
meromorphic
nowhere vanishing 1-form on $\bar {\cal A}_{X}$ without poles on $
{\cal A}_{X}$.
Now $\bar \Psi : \bar X -\r \bar {\cal A}_{X}$ is a dominant rational
map, and so
the meromorphic form $\omega$ pulls back via $\bar \Psi^{*}$ to a
nonvanishing
meromorphic 1-form   $\bar \Psi^{*}\omega$ which is holomorphic over
$X$.
By construction  we have that
$f^*(\bar \Psi^{*}\omega) \equiv 0$. Let $S \subset \bar X$ be the
divisor given by the meromorphic 1-form $\bar \Psi^{*}\omega$, namely the divisor
given by the $f_i, i\in I$ with respect to an open covering 
$U_i,i\in I$ such that $f_i \cdot (\bar \Psi^{*}\omega)|_{U_{i}}$ 
is a holomorphic 1-form
with isolated singularities only on $U_i$ for all $i\in I$. If
$f(\cz) \subset S$ then $f$ is algebraically degenerate. Otherwise
$f$ is contained in a leaf of the foliation given by the 1-form 
$\bar \Psi^{*}\omega$, namely by the $f_i \cdot (\bar \Psi^{*}\omega)|_{U_{i}}$.
 By Theorem
\ref{t15} it
follows that $f$ is algebraically degenerate, which finishes the proof
of
Proposition \ref{p16}. \qed

\subsection{The case of non-dominant quasi-Albanese map}

The following is a direct consequence of the log-Bloch theorem and the
universal property of the quasi-Albanese map.

\begin{prop}\label{nond}
Let $(\bar X,D)$ be a log surface  with log
irregularity $\bar q_X=2$.  Assume that the compactified quasi-Albanese map
$\bar \alpha_X:\bar X - \r \bar {\cal A}_X$ is not dominant. Then 
$Y=\bar\alpha_X(X)\cap {\cal A}_X$ is a hyperbolic algebraic curve.
Hence every entire curve $f:\cz \r X$ is algebraically
degenerate.
\end{prop}

\noindent{\bf Proof } The image $\bar Y :=\bar \alpha_{X}
(\bar X) \subset \bar
{\cal A}_{X}$ of $\bar X$ under the compactified quasi-Albanese map
$\bar \alpha_{X}: \bar X \r \bar {\cal A}_{X}$ is a proper algebraic
subvariety.
Since $\bar q_{X}=2$, so in particular $X$ admits nontrivial log
1-forms,  $\bar Y$
cannot degenerate to a point,  and so we have $\dim \bar Y =1$. By the log
Bloch's theorem
due to Noguchi '81 \cite{No2},  the Zariski closure of the image of
an entire curve in ${\cal A}_X$ is a translate
of an algebraic subgroup of ${\cal A}_{X}$. But by the universal
property of the quasi-Albanese
map (see Iitaka '76 \cite{Ii2}), $Y=\bar Y \cap {\cal A}_{X}$ cannot be
such a translate.
Hence, $Y$ is a hyperbolic algebraic curve.
So the map
$\bar \alpha_{X} \circ f$ is constant, and, hence, $f$ is algebraically
degenerate.\qed

   \bigskip

\section{Proof of Theorem \ref{t2} and some applications}

\subsection{Reduction of the proof}

\noindent
We will  reduce the proof of  Theorem \ref{t2} by a result of McQuillen and ElGoul 
(which is Theorem~\ref{t15} above)
and by log-Bloch's theorem to the following.
\vspace{0.3cm}

\noindent
{\bf Claim}: {\it  Let $(\bar X,D)$ be a log surface  with $\bar q_X=2$ 
and dominant $\bar \alpha_X$ and with log Kodaira dimension $\bar
\kappa_{X}=2$. Let $f : \cz \r X$ be a Brody curve. 
Assume that $f$ is not algebraically degenerate. Then 
$\alpha_X\circ f: \cz \r {\cal A}_X$ is a translate of a complex
one parameter subgroup of ${\cal A}_X$.}
\vspace{0.3cm}  

For if $\bar\alpha_X$ is not dominant, then Theorem~\ref{t2} follows from 
Proposition~\ref{nond}. If $\bar\alpha_X$ is dominant and $f$ is not algebraically
degenerate, then by the Claim
$\alpha_X\circ f: \cz\r  {\cal A}_X$ is a translate of a complex 1-parameter subgroup
of ${\cal A}_X$. Then  Theorem~\ref{t2} follows from Proposition~\ref{p16}.

So it suffices to prove the Claim for the various cases of $q_{\bar X}$ below.

\subsection{The case $q_{\bar X}=2$}

    The euclidean metric of the universal
cover $\cz^2$ of the Albanese torus ${\cal A}_{\bar X}$ descends to a
metric $h$ on it.
We may choose a hermitean metric $g$ on $\bar{X}$ such that
$\alpha_{\bar X}^* h\le g$ and we have
$$\cz \stackrel{f}{\r} (\bar X , g) \stackrel{\alpha_{\bar X}}{\r}
({\cal A}_{\bar X} , h) \leftarrow (\cz^2, {\rm eucl.}).$$
Now since $f$ is a Brody curve, we have $|f'|_g \leq C$. By
composing with the Albanese map,
we therefore get
$$\big| (\alpha_{\bar X} \circ f)'\big |_h \leq C.$$
After lifting to $\cz^2$, the components of $(\alpha_{\bar X} \circ
f)'$ are bounded
holomorphic functions. Hence, by Liouville's theorem, they are
constant.  So
    the map
$\alpha_{\bar X} \circ f:\cz \r {\cal A}_{X}$ is a translated complex 1-parameter subgroup
of ${\cal A}_X$.
  \qed

\subsection{The case $q_{\bar X}=1$}

We have the following diagram (see (5)):

$$\begin{array}{ccccc}
\cz
& \stackrel{f}{ \r }
& \bar X
& \stackrel{\bar \alpha_X}{- \r }
& \bar {\cal A}_X \\
& & & & \\
&
&
& \stackrel{ \alpha_{\bar X}}{ \searrow }
& \downarrow \\
& & & & \\
&
&
&
& {\cal A}_{\bar X}
    \end{array}
$$
As in the case $q_{\bar X}=2$, we get that $\alpha_{\bar X} \circ f$ is
linear with respect to the coordinates from the universal
covering $\cz \r {\cal A}_{\bar X}$. If $\alpha_{\bar X}\circ f$ is constant,
then $f$ is algebraically degenerate and we are done. So we assume
from now on that this linear function is nonconstant.

Let $I \s D \s \bar X$ be the (finite) set of points of indeterminacy
of $\bar \alpha_X$,
$U \s {\cal A}_{\bar X}$ a neighborhood of the finite set $\alpha_{\bar X}(I)$
consisting of small disks around each point of $\alpha_{\bar X}(I)$.
Let $V=\alpha_{\bar X}^{-1}(U)$
and $W=f^{-1}(V)$.
Since $\bar \alpha_X$ is a morphism on the compact set $\bar X
\setminus V$,
     $(\bar \alpha_X \circ f)'$ is bounded
on $\cz \setminus W$ with respect to any hermitian metric $h$ on $\bar
{\cal A}_{X}$.

Brody curves are of order $\leq 2$. Hence, by Lemma \ref{l13}, 
$$\alpha_{X} \circ f : \cz \r {\cal A}_{X} \s \bar {\cal A}_{X}$$ is of
order $\leq 2$.

Let $\cz \times \pr_{1} \r \bar {\cal A}_{X}$ be the universal cover,
and let
$$\widetilde{(\alpha_{X} \circ f)}: \cz \r \cz \times \cz^{*} \subset
\cz \times \pr_{1}$$
be a lift of the map $\alpha_{X} \circ f : \cz \r {\cal A}_{X} \s
\bar {\cal A}_{X}$
to $\cz \times \pr_{1}$. If ${\rm pr}_{1}: \cz \times \pr_{1} \r \cz$
respectively
${\rm pr}_{2}: \cz \times \pr_{1} \r \pr_{1}$ denote the projections to
the first respectively
the second factor, we get that ${\rm pr}_{1} \circ
\widetilde{(\alpha_{X} \circ f)}$
is a lift of
$\alpha_{\bar X} \circ f$ through the universal cover $\cz \r {\cal
A}_{\bar X}$, which we know already to be a linear function. Define
$\Phi := {\rm pr}_{2} \circ \widetilde{(\alpha_{X} \circ f)}: \cz \r
\cz^{*} \subset \pr_{1}$.
If we prove that $\Phi$ is of the form $\Phi (z) = \exp (P(z))$ with a
linear polynomial $P(z)$, then
$\alpha_{\bar X} \circ f:\cz \r {\cal A}_{X}$ is a translated complex 1-parameter subgroup
of ${\cal A}_X$.

By Proposition \ref{p9},
there exists a metric $h$ on $\bar {\cal A}_X$ such that the universal
covering map $(\cz \times \pr_1, {\rm eucl.} \times {\rm FS})
\r (\bar {\cal A}_X, h)$ is a local isometry.  By this, we get the
existence of  a constant
$C$ such that
\begin{equation} \label{n1}
|(\Phi|_{\cz \setminus W})'|_{FS} \leq C
\end{equation}
Furthermore, we get the following estimate for the characteristic
function.
$$ T_{\Phi}(r, \omega_{{\rm FS}}) =
\int_{1}^{r}\frac{dt}{t}\int_{|z|<t}\Phi^{*}\omega_{{\rm FS}}
\leq \int_{1}^{r}\frac{dt}{t}\int_{|z|<t}\widetilde{(\alpha_{X}\circ
f)}^{*}
\omega_{{\rm eucl.} \times {\rm FS}}$$
$$= \int_{1}^{r}\frac{dt}{t}\int_{|z|<t}(\alpha_{X}\circ
f)^{*}\omega_{h}
=T_{(\alpha_{X} \circ f)}(r, \omega_{h})
$$
Hence, from $\rho_{(\alpha_{X} \circ f)} \leq 2$, it follows
$\rho_{\Phi} \leq 2$.
So $\Phi  (z) = \exp ({P(z)})$ with $\deg P \leq 2$. If $\deg P
\leq
1$, the proof is finished.
    So we may assume $\deg P=2$.
    By a linear transformation $z \mapsto az+b$ in $\cz$ we may assume that
$P(z)=z^{2}+c$. Then,  by a multiplicative transformation $w \mapsto
\exp (c) \cdot w$,
we may assume $$P(z)=z^{2}\,.$$
 Since $U \s {\cal A}_{\bar X}$ is a small neighborhood of the finite set
$\alpha_{\bar X}(I)$ and the map $\alpha_{\bar X} \circ f$ is {\bf
linear} with respect to
the coordinates from the universal cover of ${\cal A}_{\bar X}$,
we get that $\alpha_{\bar X} \circ f$ is also a universal covering map.
Hence, up to a translation, $\alpha_{\bar X} \circ f$ 
is a group morphism with kernel $\Gamma\s \cz$.
Then
$W=(\alpha_{\bar X} \circ f)^{-1}(U)$ is the union of the translations
by the lattice $\Gamma$ of a finite number
of small disks in $\cz$.

Hence, there exists a sequence on the diagonal
$$(z_v=x_v + ix_v)_{v \r \infty} \s \cz \setminus W\; {\rm with} \; x_v
\r \infty \; .$$
\vspace{5cm}\\
We have $$|\Phi'|_{FS}(z)= \frac{|2z|\exp ({x^2-y^2})}{1+\exp
\left( {2(x^2-y^2)}\right) }$$ and, hence, since
$\Phi'$ is bounded on $\cz \setminus W$ by (\ref{n1}):
$$C \geq |\Phi'|_{FS}(z_v)= \frac{|2z_v|\exp ({0})}{1+\exp
({0})} \r \infty$$
This is a contradiction since $|z_{v}| \r \infty$. So the assumption
deg$P=2$ was wrong. \qed

\subsection{The case $q_{\bar X}=0$}
Let $(\bar X, D)$ be the log surface and assume
$f:\cz \r X$ is a Brody curve which is
not algebraically degenerate
with $f(\infty)$ its
limit set. We recall (subsection 2.2) that for a rational function $\bar
\Phi: \bar X -
\r \pr_{1}$ which extends a  holomorphic function $\Phi: X \r \cz^{*}$,
we always have $\Phi \circ f (z)= \exp (P(z))$, with $P(z)$ a
polynomial of degree
at most 2.
\begin{lemma}\label{l17}
Assume there exists a rational function $\bar \Phi: \bar X
- \r \pr_{1}$ which extends a  holomorphic function
$\Phi: X \r \cz^{*}$ such that $\Phi \circ f (z)= \exp (P(z))$, with
$P(z)$ a
polynomial of degree (exactly) $2$. Then we have $f(\infty)=
\bigcup_{j=1}^{l}D_{j}$, where
$D_{1},...,D_{l}$ is a subset of the set of the irreducible components
$D_{1},...,D_{k}$ of $D$. Moreover, for any such rational
function $\bar \Phi$,
we have
$\bar \Phi (D_{j} \setminus I_{\bar \Phi}) \equiv 0\; {\rm or}\; \equiv
\infty$ for $j=1,...,l$,
where $I_{\bar \Phi}$ denotes the set of points of indeterminacy of
$\bar \Phi$.
\end{lemma}
\noindent {\bf Proof }
Let $\bar \Phi: \bar X -
\r \pr_{1}$ be a rational function which extends a  holomorphic function
$\Phi: X \r \cz^{*}$ such that $\Phi \circ f (z)= \exp (P(z))$, with
$P(z)$ a
polynomial of degree (exactly) 2, and let
$I_{\bar \Phi}$ be the set of points of indeterminacy of $\bar \Phi$.
We first prove
\begin{equation}\label{14}
\bar \Phi (f(\infty) \setminus I_{\bar \Phi}) \subset \{0, \infty \}
\subset \pr_{1}\,.
\end{equation}
We may assume by a linear coordinate change in $\cz$ and a multiplicative
transformation in
$\cz^{*} \subset \pr_{1}$, that
$$P(z) =  z^2 \,.$$
Assume that (\ref{14}) does not hold. Then
there exists a point
$p \in f(\infty)\setminus I_{\bar \Phi}$ such that
$\bar \Phi (p) \in \cz^{*}$. Let $U(p) \subset \bar X$ be a
neighborhood such that
its topological closure (with respect to the usual topology) $\bar
U(p)$ does not
contain any points of $I_{\bar \Phi}$. Note that $\bar \Phi$ is a
holomorphic
function in a neighborhood of $\bar U(p) \subset \bar X$. Since $p \in
f(\infty)$,
there exists a sequence $(z_{v})=(x_{v}+iy_{v})$, $v \in \nz$, such that
$|z_{v}| \r \infty$ and $f(z_{v})\r p$. Without loss of generality, we
may assume
that $f(z_{v}) \in U(p)$ for all $v \in \nz$. Then we have
\begin{equation}\label{15}
\exp (x_{v}^{2}-y_{v}^2) = |\exp (x_{v}^{2}-y_{v}^{2}) \cdot
\exp (2ix_{v}^{2}y_{v}^{2})| = |\exp (z_{v}^{2})|
\end{equation}
$$= |\Phi \circ f
(z_{v})|
\r |\bar \Phi (p)| = C_{1} > 0$$

Since $f$ is a Brody curve, its derivative is uniformly bounded on
$\cz$, and since
    $\bar U(p) \subset \bar X$ is compact, the derivative of
    $\bar \Phi|_{\bar U(p)}$ is bounded too. Hence, there exists a
constant $C_{2}>0$
    such that $|(\Phi \circ f)'|_{FS} \leq C_{2}$ on $f^{-1}(\bar U(p))$.
    So in particular
    \begin{equation}\label{16}
    C_{2}\geq |(\Phi \circ f)'(z_{v})|_{FS} =
    \frac{|2z_{v}| \cdot |\exp (z_{v}^{2})|}{1+|\exp (z_{v}^{2})|^{2}}
    = \frac{2|z_{v}| \cdot \exp (x_{v}^{2}-y_{v}^{2})}
    {1+\left( \exp (x_{v}^{2}-y_{v}^{2})\right)^{2}} \r \infty
    \end{equation}
    by (\ref{15}) and as $|z_{v}| \r \infty$. This contradiction proves
(\ref{14}).

    Now assume that there exists a rational function $\bar \Phi: \bar X -
\r \pr_{1}$ which extends a  holomorphic function
$\Phi: X \r \cz^{*}$ such that $\Phi \circ f (z)= \exp (P(z))$, with
$P(z)$ a
polynomial of degree (exactly) 2. By (\ref{14}) we have
$$\left( f(\infty) \setminus I_{\bar \phi}\right) \subset \bar
\Phi^{-1}(\{ 0,
\infty \})
\subset D \subset \bar X.$$
Now the Lemma follows from Proposition \ref{p14}.
\qed
\vspace{3mm}

For the rest of this subsection, we assume $\bar q_{X}=2$ and $q_{\bar
X}=0$.
Let $\bar \Phi_{i}=( \bar \alpha_{X})_{i}: \bar X - \r \pr_{1}$,
$i=1,2$,
be the two components of the quasi-Albanese map (see subsection 2.1).
\begin{lemma}\label{l18}
Let $M= \left(
\begin{array}{cc}
m_{11} &m_{12}\\
m_{21} & m_{22}
\end{array} \right)
\in M(2\times 2, \gz)$ be a nonsingular matrix (not necessarily
invertible
over $\gz$).
Then the map $$(\bar \Psi_{1}, \bar \Psi_{2})=
( \bar \Phi_{1}^{m_{11}} \bar \Phi_{2}^{m_{12}},
    \bar \Phi_{1}^{m_{21}} \bar \Phi_{2}^{m_{22}}):\bar X - \r
(\pr_{1})^{2} $$
    is a dominant rational map extending the dominant morphism
\begin{equation}\label{17}(\Psi_{1},  \Psi_{2}):  X  \r (\cz^{*})^{2};
x \mapsto  \left( \exp \left( \int_{x_{0}}^{x}m_{11}\omega_{1}+m_{12}
\omega_{2}\right),
\exp \left( \int_{x_{0}}^{x}m_{21}\omega_{1}+m_{22}
\omega_{2}\right)\right).
\end{equation}
\end{lemma}
(We remark that the map $(\Psi_{1}, \Psi_{2})$ is not the quasi-Albanese
map in general, but
it factors through the quasi-Albanese map
by a finite \'etale map.)\\

\noindent {\bf Proof } It suffices to prove that the morphism
$(\Psi_{1},  \Psi_{2}):  X  \r (\cz^{*})^{2}$ is dominant, the rest is
clear from
the properties of the quasi-Albanese map.
Let $\omega_{1}, \omega_{2}$ be the two linearly
independent log forms corresponding to $\Phi_{1}, \Phi_{2}$.
By construction $(\Psi_{1},  \Psi_{2}):  X  \r (\cz^{*})^{2}$ is of the
form
\begin{equation}\label{18}
(\Psi_{1} , \Psi_{2})(x)= \left(\exp \left(\int_{x_{0}}^{x}m_{11}
\omega_{1}+m_{12}\omega_{2}\right),
\exp \left(\int_{x_{0}}^{x}m_{21}\omega_{1}+m_{22}
\omega_{2}\right)\right)
\end{equation}
where $x_{0}$ is a fixed point in $X$. It suffices to prove that the
lift of this map
to the universal cover,
$$X \r \cz^{2}; \; x \mapsto
\left(m_{11}\left(\int_{x_{0}}^{x}\omega_{1}\right)+m_{12}\left(
\int_{x_{0}}^{x}\omega_{2}\right),
m_{21}\left(\int_{x_{0}}^{x}\omega_{1}\right)+m_{22}\left(
\int_{x_{0}}^{x}\omega_{2}\right)\right)
$$ has rank two in some points. But this is true since $M$ is non
singular and since the lift
to the universal cover of the quasi-Albanese map, this means
$$X \r \cz^{2}; \; x \mapsto
\left(\int_{x_{0}}^{x}\omega_{1}, \int_{x_{0}}^{x}\omega_{2}\right)$$
has this property.
\qed
\vspace{3mm}

By Lemma \ref{l13}, we have $\Psi_{i} \circ f = \exp
\left(Q_{i}(z)\right)$ with
$\deg Q_{i} \leq 2$.
Then if one of the $Q_{i}$ has degree $0$, it follows immediately that
$f$
is algebraically degenerate (since $\Psi_{i} \circ f $ is constant
then).
Our aim is now to choose the nonsingular matrix $M$ such that $\deg
Q_{i} \leq 1$
for $i=1,2$. Then if one of the $Q_{i}$ has degree $0$, one is done.
If both polynomials
$Q_{i}$ have degree $1$,
then $(\Psi_1,\Psi_2)\circ f : \cz \r \cz^*\times\cz^*$
is a translate of a complex 1-parameter subgroup of $\cz^*\times\cz^*$.
But we have the commutative diagram

$$\begin{array}{ccc}
 \cz^2
& \stackrel{M}{ \r }
& \cz^2 \\
 & & \\
\downarrow &
& \downarrow \\
& & \\
{\cal A}_X& \rightarrow
& (\cz^*)^2
    \end{array}
$$
where $M$ is a linear isomorphism given by the matrix $M$ and the
vertical arrows are the universal covering maps. So
$\alpha_{\bar X} \circ f:\cz \r {\cal A}_{X}$ is also a translated complex 1-parameter subgroup
of ${\cal A}_X$.

By Lemma \ref{l13} we have
$\Phi_{i} \circ f = \exp (P_{i}(z))$ with $\deg P_{i} \leq 2$ for
the two components
of the quasi-Albanese map. If  $\deg P_{i} \leq 1$ for $i=1,2$, we
put $M=I_{2}$,
the identity matrix, and the proof is finished. Otherwise we may assume
without loss of generality that
    $\deg  P_{2} = 2$. Then by Lemma \ref{l17}, we have  $f(\infty)=
\bigcup_{j=1}^{l}D_{j}$.
    Consider the submatrix of the residues of $\omega_{1}, \omega_{2}$
with respect to the first $l$
    divisors forming $f(\infty)$
    \begin{equation}\label{19}
\begin{array}{ccccc}
&
D_1&
D_2 &
... &
D_l \\
\omega_1 &
a_{11} &
a_{12}  &
... &
a_{1l} \\
\omega_2 &
a_{21} &
a_{22}  &
... &
a_{2l}
\end{array}
\end{equation}
    If this matrix has rank $2$, there exists a nonsingular matrix $M$
with coefficients in $\gz$
    such that after passing from the log forms $\omega_{1}, \omega_{2}$
corresponding
    to the map $(\Phi_{1}, \Phi_{2})$, to the forms
    $m_{11}\omega_{1}+  m_{12} \omega_{2}, m_{21}\omega_{1}+  m_{22}
\omega_{2} $
    corresponding to $(\Psi_{1}, \Psi_{2})$, the matrix of residues
(\ref{19}) has at least one residue
    $=0$ in every line.

    If the matrix of residues (\ref{19}) of $(\Phi_{1}, \Phi_{2})$ has
rank $1$, we can choose
    the nonsingular matrix $M$ with $m_{11}\not= 0$, $m_{21}=0$ and
$m_{22}=1$
     such that after passing from the log forms $\omega_{1}, \omega_{2}$
corresponding
    to the map $(\Phi_{1}, \Phi_{2})$, to the forms
    $m_{11}\omega_{1}+  m_{12} \omega_{2}, m_{21}\omega_{1}+  m_{22}
\omega_{2} $
    corresponding to $(\Psi_{1}, \Psi_{2})$ the matrix of residues
(\ref{19}) has all residues
    $=0$ in the first line. Here, there is the difficult case where the
resulting matrix does not have any zero in the second line.

    If the matrix of residues (\ref{19}) of $(\Phi_{1}, \Phi_{2})$ has
rank $0$,
    we just take $M=I_{2}$ the identity matrix.

    In all cases, except the difficult one, there exist some residues
    $=0$ in every line of the residue matrix for $(\Psi_{1}, \Psi_{2})$.
     So by Lemma \ref{l17} and equation (\ref{8}),  we have $\deg Q_{i}
\leq 1$
for $i=1,2$, and we are done.

So we are left with the only (difficult) case that the matrix of
residues of
$(\Psi_{1}, \Psi_{2})$ with respect to the log forms
$\tilde \omega_{1}:= m_{11}\omega_{1}+m_{12}\omega_{2}$ and $\omega_{2}$
looks as follows (with $a_{2j}\not= 0$ for all $j=1,...,l$)
    \begin{equation}\label{20}
\begin{array}{ccccc}
&
D_1&
D_2 &
... &
D_l \\
\tilde \omega_1 &
0&
0  &
... &
0 \\
\omega_2 &
a_{21} &
a_{22}  &
... &
a_{2l}
\end{array}
\end{equation}
and that  $\Psi_{i} \circ f = \exp \left(Q_{i}(z)\right)$ with $\deg
Q_{1}
\leq 1$ (which follows again by Lemma \ref{l17} and equation (\ref{8}))
and $\deg Q_{2} = 2$. In this case, we have to use the explicit
geometry
of the entire curve $f$ with respect to the map $(\Psi_{1}, \Psi_{2})$
in a similar way as we did in
the proof for the case $q_{\bar X}=1$ above.

\begin{lemma}\label{l19}$\ $\\
a) The rational function $\bar \Psi_{1}: \bar X - \r \pr_{1}$ is a
morphism in a
neighborhood of $f(\infty)= \bigcup_{j=1}^{l}D_{j} \subset \bar
X$.\\[1mm]
b) There exists at most one irreducible component  of $f(\infty)$, say
$D_1$,
such that
   $\bar \Psi_{1}(D_{1})=\pr_1$. If it does not exist, then 
for all points $x\in D_i$   we have $\bar \Psi_{1}(x) \in \cz^{*}$.
   If it exists, there exists exactly one point $x_{0}\in f(\infty)$ and
exactly one point
$x_{\infty} \in f(\infty)$, both lying in $D_1$, such that $\bar \Psi_{1}(x_{0})=0$ and
$\bar \Psi_{1}(x_{\infty})=\infty$ (in particular, for all points
$x \in f(\infty)\setminus \{x_{0}, x_{\infty}\}$ , we have $\bar \Psi_{1}(x) \in \cz^{*}$).
Furthermore, $x_{0}$ respectively $x_{\infty}$ are intersection points
of $D_{1}$ with components
$D_{j}$ with $j \geq l+1$ (meaning components not belonging to
$f(\infty)$) with residues
$a_{1j}>0$ respectively $a_{1j}<0$.
\end{lemma}
\noindent {\bf Proof }
For a): This follows immediately from  Proposition \ref{p12}
since by (\ref{20}), all the residues of $\bar \Psi_1$ on
all components of
$f(\infty)$ are $=0$.\\
For b): Since $\bar \Psi_{1}$ is a morphism around $f(\infty)$, the
image of every
irreducible component of  $f(\infty)$ is a (closed) subvariety in $\pr_{1}$, and
hence a single point or
all of $\pr_{1}$. Since the residue of $\bar \Psi_1$ on each component of
$f(\infty)$ is $0$,  it follows by
(\ref{8}) that the image of each component, if it is a single point, has
to lie in
$\cz^{*}$. By the same token, we also get
that any point $x \in f(\infty)$ which is mapped to $0$ respectively to
$\infty$ by $\bar \Psi_{1}$
has to be an intersection point with a component of $D$ having residue
$>0$ respectively $<0$,
and which  cannot be one of the components of $f(\infty)$. Since $D$ is
a normal
crossing divisor, no three irreducible components meet in one point,
and so
no  point mapping to $0$ or to $\infty$ can lie on the intersection of
different irreducible
components of $f(\infty)$. We will now prove that
there can exist at most one point $x=x_{0} \in f(\infty)$ with $\bar
\Psi_{1}(x)=0$
and at most one point $x=x_{\infty} \in f(\infty)$ with $\bar
\Psi_{1}(x)=\infty$.
   From this it follows immediately that there can be at most one
component of
$f(\infty)$
mapping onto $\pr_{1}$, which then has to contain both $x_0$ and $x_{\infty}$.

Let $\{(x_{0})_{1},...,(x_{0})_{m}\}= \{x \in f(\infty): \bar
\Psi_{1}(x)=0\}$.
Take a neighborhood $U(f(\infty))$ such that $\bar \Psi_{1}$ is still a
morphism in a neighborhood
of the closure $\bar U(f(\infty))$ of $U(f(\infty))$ (this is possible
since there is only a
finite number of points of indeterminacy  of $\bar \Psi_{1}$). Take
small neighborhoods
$U_{p}=U_{p}\left((x_{0})_{p}\right)$, $p=1,...,m$, such that their
closures are
still contained in
$U(f(\infty))$.

Assume $m \geq 2$. For $p=1,2$, take sequences $(z_{v}^{(p)})$, $v \in
\nz$, such that
$|z_{v}^{(p)}| \r \infty$ and $f(z_{v}^{(p)}) \r (x_{0})_{p}$. Without
loss of generality,
we may assume that $f(z_{v}^{(p)}) \in U_{p}$ for all $v \in \nz$.
Let $[z_{v}^{(1)}, z_{v}^{2}] \subset \cz$ be the (linear) segment
between
$z_{v}^{(1)}$ and $z_{v}^{(2)}$ in $\cz$. Then there exists a point
$$w_{v} \in ([z_{v}^{(1)}, z_{v}^{(2)}] \cap f^{-1}(U(f(\infty))))
\setminus
\bigcup_{p=1}^{m}f^{-1}(U_{p})\, :$$ The image of the segment
$[z_{v}^{(1)}, z_{v}^{(2)}]$
under $f$ joins the two points $f(z_{v}^{(1)}) \in U_{1}$ and
$f(z_{v}^{(2)}) \in U_{2}$.
Since the two neighborhoods are disjoint, it has to cross the boundaries
of both of them.
Let $w_{v}$ be
the inverse image of such a crossing point with one of these boundaries.

By (\ref{20}) we have $\bar \Psi_{1} \circ f (z) = \exp (az+b)$ with
$a,b \in \cz$.
Since $f(z_{v}^{(p)}) \r (x_{0})_{p}$ we get
$\bar \Psi_{1} \circ f ( z_{v}^{(p)}) \r \bar \Psi_{1}((x_{0})_{p})=0$,
$p=1,2$.
So $\exp (az_{v}^{(p)}+b) \r 0$, meaning ${\rm Re} (az_{v}^{(p)}+b) \r
- \infty$, $p=1,2$.
Now $w_{v} \in [z_{v}^{(1)}, z_{v}^{(2)}]$, so there exist $\lambda_{v}$
with
$0 \leq \lambda_{v} \leq 1$ such that
$w_{v}=\lambda_{v}z_{v}^{(1)}+(1-\lambda_{v})
z_{v}^{(2)}$. Then we have
\begin{equation}\label{21}
{\rm Re} (aw_{v}+b) = {\rm Re} (\lambda_{v}
(az_{v}^{(1)}+b)+(1-\lambda_{v})(az_{v}^{(2)}+b))
\end{equation}
$$= \lambda_{v} {\rm Re} (az_{v}^{(1)}+b) +(1-\lambda_{v}) {\rm Re}
(az_{v}^{(2)}+b) \r - \infty$$
since $\lambda_{v}, 1-\lambda_{v}\geq 0$. In particular $|w_{v}| \r
\infty$.

Consider the sequence $f(w_{v})$, $v \in \nz$. After passing to a
subsequence, we get
$f(w_{v}) \r P \in \bar X$. By construction, we have $P \in
f(\infty)$,
and by (\ref{21}) we get
\begin{equation}\label{22}
\bar \Psi_{1}(P)= \lim_{v \r \infty} \bar \Psi_{1}(f(w_{v})) =
\lim_{v \r \infty} \exp (aw_{v}+b) =0\,.
\end{equation}
But $w_{v} \in f^{-1}(U(f(\infty))) \setminus
\bigcup_{p=1}^{m}f^{-1}(U_{p}((x_{0})_{p}))$,
so $P\not= (x_{0})_{v}$, $v=1,...,m$. This is a contradiction, since
$\{(x_{0})_{1},...,(x_{0})_{m}\}= \{x \in f(\infty): \bar
\Psi_{1}(x)=0\}\ni P$.
So our assumption $m \geq 2$ was wrong, and there exists at most one
point
$x_{0}\in f(\infty) $ such that $\bar \Psi_{1}(x_{0})=0$.

The proof for $x_{\infty}$ is exactly the same, just change $0$ and
$\infty$ in the
proof above, and in (\ref{21}), change $-\infty$ to $+\infty$.
\qed

Now we can finish up the proof of the claim:
If there does not exist an irreducible component $D_1$ of
$f(\infty)$ such that $\bar \Psi_1(D_1)=\pr_1$ (case 1),
we leave $(\Psi_1,\Psi_2)$ as it is. If it does exist (case 2),
then we may assume 
 $x_{0} \in D_{1} \cap
D_{l+1}$. The matrix of residues
(\ref{20}) extends to
    \begin{equation}\label{20a}
\begin{array}{cccccc}
&
D_1&
D_2 &
... &
D_l &
D_{l+1}\\
\tilde \omega_1 &
0&
0  &
... &
0&
a_{1,l+1} \\
\omega_2 &
a_{21} &
a_{22}  &
... &
a_{2l}&
a_{2,l+1}
\end{array}
\end{equation}
with $a_{1,l+1}>0$. We put
$$ \tilde \omega_{2}:=a_{1,l+1}\omega_{2} - a_{2,l+1}\omega_{1}$$
Then the new matrix of residues becomes
    \begin{equation}\label{20b}
\begin{array}{cccccc}
&
D_1&
D_2 &
... &
D_l &
D_{l+1}\\
\tilde \omega_1 &
0&
0  &
... &
0&
a_{1,l+1} \\
\tilde \omega_2 &
a_{21} &
a_{22}  &
... &
a_{2l}&
0
\end{array}
\end{equation}

By Proposition \ref{p12}, we get that $x_{0}\in D_{1}\cap D_{l+1}$
cannot be a
point of indeterminacy of $\bar \Psi_{2}$, since the corresponding
product of residues
is zero.
In both cases we get, by a linear coordinate change in $\cz$ and a multiplicative
transformation in
     $\cz^{*}\subset \pr_{1}$, that
$$(\Psi_1,\Psi_{2})\circ f(z) = (\exp (az),\exp (z^2))\, , \, \ \  a\in \cz^*.$$
By Lemma~\ref{limit1} and Lemma~\ref{limit2}, we get 
\begin{equation}\label{eq*}
\{0\}\times \pr_1\s (\Psi_1,\Psi_2)\circ f(\infty)\s (\bar\Psi_1,\bar\Psi_2)(f(\infty)).
\end{equation}
We will reach a contradiction to (\ref{eq*}) in both cases. In case 1, we have
$\bar\Psi_1(f(\infty))\s \cz^*$, so $$(\bar\Psi_1,\bar\Psi_2)(f(\infty))\cap (\{0\}\times \pr_1)
=\emptyset,$$ contradicting (\ref{eq*}). In case 2, we get that $(\bar\Psi_1|_{f(\infty)})^{-1}(0)=\{x_0\}$
and $x_0$ is not a point of indeterminacy  of $\bar\Psi_2$. So we have 
$$(\bar\Psi_1,\bar\Psi_2)(f(\infty))\cap (\{0\}\times \pr_1)= 
\{ 0\}\times \{\bar\Psi_2(x_0)\}$$
which is strictly contained in $\{0\}\times \pr_1$, contradicting (\ref{eq*}).

This concludes the proof of the claim, and, by the reduction in subsection 3.1, the proof of Theorem~\ref{t2}.
\qed

\subsection{Some applications}
Theorem~\ref{t2} gives us the following corollaries on hyperbolicity.\\[-6mm]
\begin{cor}\label{c3}
    Let $(\bar X,D)$ a log surface  with log
irregularity $\bar q_X=2$ and log Kodaira dimension $\bar \kappa_{X} > 0$. 
Suppose that $X$ does not contain non-hyperbolic algebraic curves and that $D$
is hyperbolically
stratified (this means that every irreducible component  minus all the
others is a hyperbolic
curve). Then $X$ is complete hyperbolic and hyperbolically imbedded in
$\bar X$.
\end{cor}
\noindent {\bf Proof } $\bar \kappa_{X} \neq 1$, otherwise the Iitaka fibration of $X$ 
would contain non-hyperbolic fibers and so $X$ would have non-hyperbolic algebraic
curves. Hence this follows from Theorem~\ref{t2}
 and from the
main result of Green '77
\cite{Gr}.
\qed \medskip
\noindent Another consequence of Theorem \ref{t2} is the "best
possible" result for algebraic
degeneracy in the three component case of complements of plane
curves
(see Dethloff-Schumacher-Wong '95 \cite{DSW1} and \cite{DSW2},
Bertheloot-Duval '01
\cite{BD}):
\begin{cor}\label{c4}
    Let $D \subset \pr_{2}$ be a normal crossing curve of degree at least four consisting
of three
    components.
     Then every Brody curve $f:\cz \r \pr_{2}\setminus D$
is algebraically
degenerate.
\end{cor}
\noindent {\bf Proof }  $\bar K_{(\pr_{2}\setminus D)} =
{\cal O}\left({(\rm deg}D)-3\right)$ is very ample and, hence,  big. Let
$D_{i}=\{f_{i}=0\}$, where
$f_i$ is a homogeneous equation for $D_i$, for
for $i=1,2,3$. Then
$$\omega_{j}=d \log \left( \frac{f_{j}^{{\rm deg}f_{3}}}
{f_{3}^{{\rm deg}f_{j}}}\right) \;,\;  j=1,2,$$
are linearly independent log forms, and it is easy to see that there
are no others.
Hence, $\bar q_{(\pr_{2}\setminus D)}=2$. \qed

 \section{Proof of Theorem~\ref{neut} and applications}

\subsection{Kawamata's theorem}

Recall that for a fibered variety, i.e.,  a dominant algebraic morphism
between irreducible, reduced quasiprojective varieties with irreducible
general fibers, we have the following addition theorem of 
Kawamata '77 (\cite{add}) in dimension two:

\begin{theo} \label{t10b} Let $f:V\r B$ be a fibered variety for a nonsingular
algebraic surface $V$ and a nonsingular algebraic curve $B$. Let $F$ be a general
fiber of $f$. Then $$ \bar\kappa_V\geq \bar \kappa_F+\bar \kappa_B.$$

\end{theo}

\noindent
Recall also the following result of Kawamata '81 \cite[Theorems 26, 27]{Kaw}:
\begin{theo} \label{t10a}
Let $X$ be a normal algebraic variety, ${\cal A}$  a semi-abelian
variety and let $f:X \r {\cal A}$
be a finite morphism. Then $\bar \kappa_{X} \geq 0$ and there exist a
connected complex  algebraic subgroup
${\cal B}\subset {\cal A}$, \'etale covers $\tilde X$ of $X$ and 
$\tilde {\cal B}$ of ${\cal B}$, and a normal algebraic variety $\tilde Y$
such that
\begin{itemize}
\item[(i)] $\tilde Y$ is finite over ${\cal A}/{\cal B}$.
\item[(ii)] $\tilde X$ is a fiber bundle over $\tilde Y$ with fiber
$\tilde {\cal B}$ and translations
by $\tilde {\cal B}$ as structure group.
\item[(iii)] $\bar \kappa_{\tilde Y} = \dim \tilde Y= \bar
\kappa_{X}$.
\end{itemize}
If further $\bar \kappa_{X} = 0$ and $f$ is surjective, then $f$ is an \'etale morphism.

\end{theo}

For the proof of Theorem~\ref{neut} we will need the following  extension of this
theorem in the case of surfaces. 

\begin{theo} \label{neutt}
Let $X$ be a normal algebraic surface, ${\cal A}$  a semi-abelian
surface and let $f:X \r {\cal A}$
be a finite morphism. 
Let $X_{0} \subset X$ be an open algebraic subvariety.

\begin{itemize}

\item[\rm{(1)}] In the case $\bar \kappa_{X_{0}}=1$, let 
$\Phi : X^{*}\rightarrow Y^{*}$ be the logarithmic Iitaka fibration of
$X_{0}$, $\Psi : X^{*} \rightarrow X_{0}$ the birational morphism
relating $X_{0}$ to $X^*$, and for $y \in Y^{*}$,
let $X_{y}^{*}=\Phi^{-1}(y)$, $X_{y}=\Psi (X_{y}^{*}) \subset X_{0} 
\subset X$ 
and $B_{y}= f(X_{y})$.  Then for generic $y \in Y^{*}$, $X_{y} \subset X$
is a closed subvariety  and $B_{y}$ is a translate of a complex one parameter
algebraic subgroup ${\cal B}$ of $\cal A$. Moreover, there are  \'etale covers
$\tilde X$ of $X$ and 
$\tilde {\cal B}$ of ${\cal B}$, and a smooth algebraic curve $\tilde Y$
such that
\begin{itemize}
\item[(i)] $\tilde Y$ is finite over ${\cal A}/{\cal B}$.
\item[(ii)] $\tilde X$ is a fiber bundle over $\tilde Y$ with fiber
$\tilde {\cal B}$ and translations
by $\tilde {\cal B}$ as structure group. In particular, $X$ and $\tilde X$ are smooth. 
\item[(iii)] $ \dim \tilde Y= \bar \kappa_{X_0}=1$.

\end{itemize}
In particular,
 for generic $y \in Y^{*},$ $X_{y} \subset X$ is equal to the image of a suitable fiber of the fiber
bundle $\tilde X$ over $\tilde Y,$ and the image of a generic fiber 
of this fiber bundle is of the form $X_{y}$ for $y \in Y^{*}$.

\item[\rm{(2)}] If $\bar \kappa_{X_{0}}=0$, then $f$ is an \'etale morphism and 
$X\setminus X_0$ is finite.
\end{itemize}

\end{theo}

  \noindent {\bf Proof } 
Since the proof uses among others  the same ideas
as that given by Kawamata in \cite{Kaw} (see Theorems 13, 23, 26 and 27 of Kawamata's
paper \cite{Kaw}), we just indicate the differences and what has to be added:

\begin{itemize}
\item[(1)] In the case $\bar \kappa_{X_{0}}=1$,
the key point is to observe that despite the composed map
$X^{*} \stackrel{\Psi}{\rightarrow} X_{0} \hookrightarrow X$ being not
proper in general, for generic $y \in Y^{*}$, $X_{y} \subset X$
is a closed subvariety. In order to see this, let $\bar X_{y} \subset \bar X$
be the (Zariski) closure of $X_{y}$ in $\bar X$ (it is the curve 
$\bar \Psi (\bar X_{y}^{*})$). Then $\bar X_{y} \cap X$ is the closure
of $X_{y}$ in $X$, $\bar B_{y}= f(\bar X_{y} \cap X)$ is the closure of 
$B_{y}\subset {\cal A}$ and we have for the logarithmic Kodaira dimensions 
(for generic $y \in Y^{*}$):
$$0 = \bar \kappa_{X_{y}^{*}} =  \bar \kappa_{X_{y}} \geq
 \bar \kappa_{\bar X_{y} \cap X}\geq  \bar \kappa_{\bar B_{y}} \geq 0\, ,$$
 where the last inequality follows since we can restrict the constant 1-forms of
 ${\cal A}$ to $\bar B_{y}$, and all the others follow by lifting back log
 multi-canonical forms through logarithmic (!) morphisms. So we get in particular
 for the inclusion of algebraic curves $X_{y} \subset \bar X_{y} \cap X$ 
 that $0 =   \bar \kappa_{X_{y}}=
 \bar \kappa_{\bar X_{y} \cap X}$ (remark that $X_{y} \subset X$ is
 algebraic since $X_{y}\subset X_{0}$ is algebraic and
  $X_{0}\subset X$ is an open algebraic subvariety). But an
 algebraic curve having log Kodaira dimension zero is either isomorphic to
 an elliptic curve or to $\cz^{*}\subset \pr_{1}$, and if we delete any
 other point of any of these algebraic curves, they become hyperbolic, so their
 log Kodaira dimension becomes one. This means 
 $X_{y} = \bar X_{y} \cap X$ is a closed subvariety of $X$. The rest of the proof
 is exactly analogous to the proof of Kawamata's Theorem 27. We observe here
 that the normal curve $\tilde Y$, having singularities only in codimension two, 
  is automatically smooth. 
  
  \item[(2)] As $f$ is a finite morphism between algebraic surfaces, it is surjective.
 So the first part of Theorem~\ref{t10a} gives  $0=\bk_{X_0}\geq \bk_X\geq 0.$
Hence the first part of our claim follows from the last part of Theorem~\ref{t10a}.
Since $X$ is finite and \'etale over $\cal A$, $X$ is also a semi-abelien variety since
it is a quotient of affine space by a discrete subgroup and the
kernel in $X$ of the map to the abelian variety quotient of $\cal A$ 
is a finite cover of that of ${\cal A}$.
Let $\bar X$ be a compactification of $X$ as given in subsection 2.2.
If $X\setminus X_0$ is not finite, then there exist a curve $C$ in $\bar X$ 
intersecting $X$ but not $X_0$. We will reach a contradiction.

Assume
that $X$ is not a simple abelian surface. 
 We follow the proof of Theorem 26 in \cite{Kaw} to obain a contradiction.
As $X$ is semi-abelian, it is a fiber bundle
over a flat curve $B$, i.e.$\,B$ has trivial log-canonical bundle so that $B$ is 
either $\cz^*$ or an elliptic curve. Moreover, the fiber bundle has flat curves
as fibers. Now the log-Kodaira dimension becomes positive if we remove
a finite nonempty set from a flat curve, whose log-Kodaira dimension is $0$.
If $C$ intersects the generic fiber
of this fiber bundle, then we have a contradiction with the addition theorem above. 
Hence, $C$ is vertical and so $X\setminus C$ is a fiber space over a curve of positive
Kodaira dimension. This again contradicts the addition theorem. 

Hence we may assume
that $X$ is a simple abelian surface. Now it is well known that $C$ is then an ample 
divisor in $X$ so that its $L$-dimension is 2. This latter fact can be seen directly as
follows \begin{itemize} \item{} If $C^2>0$, then Riemann-Roch gives that
$$h^0(mC)-h^1(mC)+h^2(mC)= O(m^2)$$
while Serre-duality gives $h^2(mC)=h^0(-mC)=0$. Hence 
$h^0(mC)= O(m^2)$ and the $L$-dimension of $C$ is $2$. If $C^2\leq 0$,
then the algebraic subgroup $G=\{t\in X\ |\ C+t=C\}$ must contain differences of points in $C$
and so is at least one dimensional. As $X$ is simple, $G=X$ which is a contradiction. --
\end{itemize}
To compute the log-Kodaira dimension of $X\setminus C$, we blow up
$X$ successively until the total reduced transform of $C$ becomes normal crossing,
where we let $E_i$ be the exceptional divisor at the $i$-th step. Then
we know that at each step,  the resulting canonical divisor is $\sum_i E_i$ and
$C= \tilde C + \sum_i m_i E_i$ as divisors ($m_i> 0$), where we 
have suppressed the pullback symbols and $\tilde C$ is the strict
transform of $C\subset \bar X$. Hence, for the final blown up variety $\hat X$ of $X$, the 
log-canonical divisor of $(\hat X, \tilde C)$ is $K=\tilde C+\sum_i E_i$. 
But this means that a multiple of
$K$ is $C$ plus an effective divisor 
and therefore the log-Kodaira dimension 
of $\hat X\setminus \tilde C$ is $2$. This is a contradiction
since then $X_0\subset \hat X\setminus \tilde C$ 
must have the same log-Kodaira dimension. 
\qed

\end{itemize}

\subsection{Stein factorization and condition (*)}

Let $(\bar X,D)$ be a log surface  with log irregularity $\bar q_X=2$ .
Let $\bar \alpha_X :\bar X - \r \bar {\cal A}_X$ be the compactified Albanese map, $I$ its
finite set of points of indeterminacy and 
$\bar \alpha_0=\bar \alpha_X|_{\bar X\setminus I}$.
We assume in this subsection that $\bar \alpha_X$ is dominant.
Recall the following condition:
\begin{itemize}
\item[\rm{(*)}] For all $z\in {\cal A}_X$ and $\bar E$ a connected component of
the Zariski closure of $\bar \alpha_0^{-1}(z)$ with $\bar E\cap X\neq \emptyset$,
any connected component of $D$ intersecting $\bar E$ is contained in $\bar E$ 
(i.e. $\bar E$ is a connected component of $\bar E\cup D$).
\end{itemize}
Condition (*) is a natural condition in the sense that all the data can be read directly
from two linearly independent log $1$-forms, $\omega_1, \omega_2$ without referring to
$\alpha_X$ or $\bar \alpha_X$. Indeed, $\alpha_X$
is not dominant iff $\omega_1\wedge \omega_2\equiv 0$ and a curve $C$ in $\bar X$ such
that neither $\omega_1$ nor $\omega_2$ is identically $\infty$ along $C$ 
is exceptional with respect to $\bar \alpha_X$ iff $i^*\omega_1\equiv 0\equiv i^*\omega_2$
where $i:C\r \bar X$ is the inclusion map.

\begin{prop} \label{stein}
Condition {\rm (*)} is equivalent to the following condition:
If $\bar{\hat{\alpha}}_{X}:\bar{\hat{X}} \r \bar {\cal A}_{X}$
is a desingularization
of the rational map $\bar \alpha_{X}$ which  is biholomorphic on
$X$, then for the Stein factorization 
$\bar{\hat{X}} \stackrel{\bar \beta}{\r}\bar Y
\stackrel{\bar \gamma}{\r} \bar {\cal A}_{X}$ 
of the morphism
$\bar{\hat{\alpha}}_{X}$ we have that
$\bar \beta (X) \subset \bar Y$ is an open algebraic subvariety.

\end{prop}

\noindent
{\bf Proof } We first prove the following lemma.

\begin{lemma}  Let $\bar{\hat{\alpha}}_{X}:\bar{\hat{X}} \r \bar {\cal A}_{X}$
be a desingularization
of the rational map $\bar \alpha_{X}$ 
via $r:\bar{\hat{X}} \r \bar X$
which  is biholomorphic on
$X$ and
$\bar{\hat{X}} \stackrel{\bar \beta}{\r}\bar Y
\stackrel{\bar \gamma}{\r} \bar {\cal A}_{X}$  the Stein factorization 
of the morphism
$\bar{\hat{\alpha}}_{X}$. We identify $X$ and $r^{-1}(X)$ and let
$\hat D$ be the reduced total transform of $D$.
Then condition {\rm (*)} is equivalent to the following condition.
\begin{itemize}
\item[\rm{(**)}] For all $z\in {\cal A}_X$ and $\bar {\hat E}$ a connected component of
$\bar {\hat \alpha}_0^{-1}(z)$ with $\bar{\hat E}\cap X\neq \emptyset$,
any connected component of $\hat D$ intersecting $\bar{\hat E}$ is contained in $\bar {\hat E}$ 
(i.e. $\bar {\hat E}$ is a connected component of $\bar{\hat E}\cup \hat D$).
\end{itemize}
\end{lemma} 
{\bf Proof of the Lemma }\\
{\bf {{(*)}} implies {{(**)}} } Assume that (**) does not hold.
Then there exist $z\in {\cal A}_X$ and a
connected component $\bar {\hat E}$ of $\bar{\hat\alpha}^{-1}(z)$ with 
nonempty intersection with $X$. Also there is
a connected component of 
$\hat D$ having nonempty intersection with but not contained in $\bar {\hat E}$. Let 
$\sum_k \hat D_k$ be its decomposition into irreducible components. We may assume
that $\hat D_1\not\s \bar {\hat E}$ so that $\bar{\hat\alpha}_X|_{\hat D_1}\not\equiv z$
and that there is a point $\hat x\in \hat D_1\cap\bar {\hat E}$.

Now $r(\bar {\hat E}) $ is a connected subset of a connected component $\bar E$ of the 
Zariski closure of $\bar\alpha_0^{-1}(z)$. We distinguish 2 cases.

\noindent
{\bf Case 1 } $r(\hat D_1)\s D$ is a curve. Then $r(\hat D_1)$ is an irreducible component
of $D$ intersecting $\bar E$ at $r(\hat x)$, but $r(\hat D_1)\not\s \bar E$
because $\bar\alpha_X|_{r(\hat D_1)\setminus I}\not\equiv z$. This contradicts
(*).

\noindent 
{\bf Case 2 } $r(\hat D_1)\s D$ is a point. In this case $r(\hat x)\in I$ is a point where
$\bar\alpha_X$ has a point of indeterminacy. By Propositions~\ref{p12} and \ref{p12a},
$r(\hat x)\in \bar E$ is an intersection point of two irreducible components $D_i$ and
$D_j$ of $D$ such that the form $\omega$ defining $\Phi$ (only around 
$r(\hat x)$ in the case of $q_{\bar X}=1$) has strictly opposite signs along
$D_i$ and $D_j$, and $\bar\alpha_X$ maps $D_i\setminus I$ and 
$D_j\setminus I$ to $\bar {\cal A}_X\setminus {\cal A}_X$ and so, in particular,
$\bar \alpha_X|_{D_i\setminus I}\not\equiv z$ and so (*) cannot hold. \\[-2mm]

\noindent
{\bf {{(**)}} implies {{(*)}} } Assume that (*) does not hold. Then there exist 
$z\in {\cal A}_X$, a connected component $\bar E$ of the Zariski closure of 
$\bar\alpha_0^{-1}(z)$ with $\bar E\cap X\neq\emptyset$ and a connected
component $D_0$ of $D$ such that $D_0\cap \bar E\neq \emptyset$ and 
$D_0\not\s \bar E$. Let $\hat D_0$ and $\bar{\hat E}$ be the proper transforms
of $D_0$ and $\bar E$  under $r:\bar{\hat X}\r \bar X$ respectively. 
Then any connected component of 
$\bar{\hat E}$ intersects $\hat D_0$ or 
is connected to $\hat D_0$ by the exceptional divisor $F$ of $r$.
Let $\hat E_0$ be such a component. Since $D_0\not\s \bar E$, we have 
$\hat D_0\not \s \bar{\hat \alpha}_X^{-1}(z)\supset \hat E_0$. So 
there is a component of $F\cup \hat D_0$ intersecting with
but not contained in $\hat E_0$. This contradicts (**). \qed
$\ $\\
We now continue with the proof of the proposition. 

Assume (**), we prove that $\bar \beta (X)\s \bar Y$ is an open algebraic
subvariety by giving an explicit description for it: By definition of the Stein
factorization, points of $\bar Y$ corresponds to connected components
of the fibers of the morphism $\bar{\hat \alpha}_X: \bar{\hat X}\r \bar{\cal A}_X$,
and $\bar\beta:\bar{\hat X}\r \bar Y$ is the canonical surjective map which
contracts these connected components to points. Since $\bar \beta$ is a
proper birational morphism, $\bar \beta(\hat D)\s \bar Y$ is algebraic.
Since $\bar \alpha_X$ is dominant, $\bar Y$ is a surface and so 
$\bar \beta(\hat D)$ can be decomposed into a  pure $1$-dimensional
subvariety $\Gamma$ and a finite set $G=\{y_1, ..., y_s\}$ in $\bar Y$. 
Let  $y\in \bar\beta (X)$, then $y=\bar\beta(\hat x)$ for some $\hat x\in X$. 
Let $z=\bar {\hat \alpha}_X(\hat x)$. By the definition of $\bar\beta$, 
$\bar\beta^{-1}(y)$ is a connected component $\bar{\hat E}$ of 
$\bar{\hat\alpha}^{-1}(z)$. By (**), all connected components of 
$\hat D$ intersecting $\bar{\hat E}$ map to the point $y$, and the 
others map to a closed algebraic set not containing $y$. This means
that $y\not\in \Gamma$ and so $\bar\beta(X)$ is just $\bar Y\setminus \Gamma$
minus a finite number of points in $G$, which is an algebraic set.

Conversely, assume that (**) does not hold. Then there exist $z\in {\cal A}_X$,
a connected component $\bar{\hat E}$ of $\bar{\hat\alpha}^{-1}(z)$ with
$\bar{\hat E}\cap X\neq \emptyset$ and an irreducible component $\hat D_0$
of $\hat D$
with $\hat x\in \hat D_0\cap \bar{\hat E}$ and $\hat D_0\not\s \bar E$. So
$\Gamma=\bar\beta(\hat D_0)\s \bar Y$ is a $1$-dimensional subvariety
containing the point $y_0=\bar\beta(\hat x)=\bar\beta(\bar{\hat E})$. So 
$\bar\beta$ has a $1$-dimensional fiber over $y_0$ and so all other
fibers over a neighborhood $U$ of $y_0$ are single points. Hence
$U\cap\Gamma\cap\bar\beta(X)=\{y_0\}$ and $\bar\beta (X) \s \bar Y$
cannot be open.\qed

\subsection{End of the proof of Theorem~\ref{neut}}

For this  proof we adopt the notations of Theorem~\ref{t10a},
meaning we
allow arbitrary boundary divisors. This is not a problem since
we do not need the condition that $f:\cz \r X$ is Brody in this 
 theorem and since we
can always lift an entire
curve $f$ through a birational map which is biholomorphic on $X$
and the property of $f$ to be algebraically
degenerate is of course invariant
under such a birational map. 

By the same token
we also may assume that $\bar \alpha_X: \bar X \r \bar {\cal A}_X$ is
(already) a morphism.
Hence, let $\bar X \stackrel{\bar \beta}{\r}\bar Y
\stackrel{\bar \gamma}{\r} \bar {\cal A}_{X}$ be a Stein factorization
of the morphism
$\bar \alpha_{X}$ and $\beta=\bar \beta|_X$. 
Then $\bar Y$ is a normal variety which
compactifies
$Y=(\bar \gamma)^{-1}({\cal A}_{X})$, and $f$ is algebraically degenerate (in $X$) if and
only
if the map
$\beta \circ f$ is algebraically degenerate (in $Y$). 

We apply Theorem~\ref{neutt}
  to  the finite morphism $\gamma=\bar\gamma|_Y : Y \rightarrow {\cal A}_X$
   and to  $ Y_{0} :=\beta (X) \subset Y$, which is an open subvariety by
   Proposition~\ref{stein}. In order to
  simplify notations, we assume that $ X$ is the
 is the total space of the Iitaka fibration of $Y_0$ (there is no problem with this
  since the entire curve $f:\cz \rightarrow X$ lifts through a birational morphism
  over $X$). So let $\Phi : X \rightarrow Z^{*}$ be the log Iitaka fibration.
Then Theorem~\ref{neutt} says that 
 there exist a
semi-abelian curve
${\cal B}\subset {\cal A}_X$, \'etale covers $\tilde Y$ and $\tilde {\cal
B}$ of $Y$
and ${\cal B}$, respectively, and a  smooth curve $\tilde Z$
such that
\begin{itemize}
\item[(1)] $\tilde Z$ is finite over ${\cal A}_X/{\cal B}$.
\item[(2)] $\tilde Y$ is a fiber bundle over $\tilde Z$ with fiber
$\tilde {\cal B}$ and translations
by $\tilde {\cal B}$ as structure group.
\end{itemize}
Moreover,  if for $z \in Z^{*}$
we put $X_{z}=\Phi^{-1}(z)$, $Y_{z}=\beta (X_{z}) \subset Y_{0} 
\subset Y$ 
and $B_{z}= \gamma(Y_{z})$, then for generic $z \in Z^{*}$, $Y_{z} 
\subset Y$
is a closed subvariety, which is equal to the image of a suitable fiber of the fiber
bundle $p: \tilde Y \rightarrow \tilde Z$, and, vice versa, the image of a generic fiber 
of this fiber bundle is of the form $Y_{z}$ for generic $z \in Z^{*}$. 

Let $\tilde Y_{0} \subset \tilde Y$ be the inverse image of
 $Y_{0}$ in  $\tilde Y$ under the \'etale cover $\tilde Y \rightarrow Y$, 
and  let $\tilde Z_{0} \subset \tilde Z $
 be the image of $\tilde Y_{0}$ under the fiber bundle projection
 $p: \tilde Y \rightarrow \tilde Z$. Then $\tilde Z_{0} \subset \tilde Z$
 is an open algebraic subcurve: In fact, since 
 $\tilde Y_{0} \subset \tilde Y$ is an open algebraic subsurface, its complement
 can only contain a finite number of curves or isolated points, which can only 
 contain a finite number of fibers of the fiber bundle 
 $\tilde Y \rightarrow \tilde Z$, so the complement of $\tilde Z_{0}$ in
 $\tilde Z$ contains at most a finite number of points. Let still 
 $\tilde W_{0} \subset \tilde Y$ be the inverse image of $\tilde Z_{0}$
 under the fiber bundle projection $p: \tilde Y \rightarrow \tilde Z$.
 
 We claim that $\tilde W_{0} \setminus \tilde Y_{0}$ is of codimension
 at least 2 in  $\tilde W_{0}$: In fact, by Theorem~\ref{neutt} the generic fiber of the fiber bundle
 $\tilde Y \rightarrow \tilde Z$ projects to a fiber 
 $Y_{z} \subset Y_{0}$, so is contained in $\tilde Y_{0}$. We know that
 since 
 $\tilde Y_{0} \subset \tilde Y$ is an open algebraic subsurface, its complement
 can only contain a finite number of (closed) curves or isolated points. But since
 such a closed curve does not hit the generic fiber of the fiber bundle (so it is contained in
 the union of a finite number of fibers since the base $\tilde Z$ is of dimension 1), it is itself
 equal to a fiber, but then it is also in the complement of $\tilde W_{0}$.
 This proves the claim.

 We now claim that $\tilde W_{0}$ is of log Kodaira dimension 
 $\bar \kappa_{\tilde W_{0}} \geq 1$:  Since $p:\tilde Y\r \tilde Z$ is a fiber bundle 
 over the smooth curve $\tilde Z$ and
 $\tilde Y$ is \'etale over $Y$, both $Y$ and $\tilde Y$ are smooth.
 Let $X^+=\bar \beta^{-1}(Y_0)\s \bar X$.  Then $X^+$ has the same log Kodaira dimension 
  as $Y_0$ as they are properly birational. But $X^+\setminus X$ lies only on the
  fibers of $\bar\beta$ above $Y_0$ and so any pluricanonical 
  logarithmic form for $(\bar X,D)$ must
  actually be a pluricanonical form over $X^+$ by the Reimann extension theorem  
  applied to $Y_0$.   Hence 
  $1=\bar \kappa_{X}  =\bar \kappa_{X^+}=\bar \kappa_{Y_{0}} \leq
 \bar \kappa_{\tilde Y_{0}}=\bar \kappa_{\tilde W_{0}}$ where the last equality
 follows again by the Riemann extension theorem.

 Next we claim that $\bar \kappa_{\tilde Z_{0}} = 1$: Let 
 $s_{1}, s_{2} \in H^{0}(\tilde Y, m \bar K_{  \tilde W_{0}})$ be two
 linearly independent log multi-canonical sections. 
 Consider the exact sequence of locally free sheaves
 $$0\rightarrow p^*\bar K_{\tilde Z_0} \rightarrow
  \bar \Omega_{\tilde W_0} \rightarrow
\bar K_{{\tilde Y}/{\tilde Z}} \rightarrow 0\, ,$$
from which it follows that 
$\bar K_{\tilde W_0}$ is isomorphic to 
$p^*\bar K_{\tilde Z_0} \otimes \bar K_{{\tilde Y}/{\tilde Z}}$.
But
$\bar K_{{\tilde Y}/{\tilde Z}}$ is dual to
$$\ker ( \bar dp : \bar T_{\tilde W_0} \rightarrow
 \pi^*\bar T_{\tilde Z_0} )$$
and this latter has a global nowhere vanishing section $v$ (that generates
the Lie group action) and, hence, is trivial. So 
$\bar K_{\tilde W_0}$ is isomorphic to $\bar K_{\tilde Z_0} $. 
Hence $s_1,s_2$ give two linearly independent sections of
$m\bar K_{\tilde Z_0} $.

 Theorem~\ref{neut} is now immediate:
 In fact, let $f:\cz \rightarrow X$ be an entire curve. If 
 $p: \tilde Y \rightarrow \tilde Z$ denotes the fiber bundle projection,
 then the entire curve $\beta \circ f: \cz \rightarrow Y_{0}$ lifts by the \'etale
 cover $\tilde Y \rightarrow Y$ to a curve having image in 
 $\tilde Y_{0} \subset \tilde W_{0}$, projecting by $p$ to an entire
 curve having image in the curve $\tilde Z_{0}$. But since 
  $\bar \kappa_{\tilde Z_{0}} = 1$,  the curve $\tilde Z_{0}$ is hyperbolic,
  so this entire curve is constant.
  Hence, the image of the entire curve $\beta \circ f$ lies in the image in $Y$
  of a fiber of the fiber bundle $\tilde Y \rightarrow \tilde Z$, and so
  the entire curve $f:\cz \rightarrow X$ is algebraically degenerate.
 \qed \vspace{1cm}

\noindent
{\bf Remark } The counter-example in Proposition~\ref{p6} shows that,
without  condition (*),
it can happen that entire curves $f:\cz \r X$ can be Zariski-dense. In fact, condition (*) does not
hold in this counter-example: With the notations as in Proposition~\ref{p6}, 
$\bar\alpha_X=\alpha_{\bar X}:\bar X\r E\times E$ is the blowup morphism
of the two point $Q_1, Q_2$ and the exceptional fibers 
$\alpha_{\bar X}^{-1}(Q_1)$ and $\alpha_{\bar X}^{-1}(Q_2)$
intersect $D$ but do not contain it.

\subsection{Some applications} 

The following result generalizes Kawamata's theorem
for normal surfaces finite over a semi-abelian surface to surfaces with
log irregularity 2 by relating it to entire holomorphic curves.

\begin{theo}\label{neupt}
Let $(\bar X,D)$ be a log surface  with log irregularity $\bar q_X=2$ .
In the case of dominant $\bar \alpha_X$, assume condition {\rm (*)}.
Then the following are equivalent.

\begin{itemize}
\item[\rm{(1)}] There is an  entire curve $f:\cz\r X$ 
such that $f^*w\equiv 0$ for some $w\in H^0(\bar T^*_X)$ and $f$ is not algebraically 
degenerate.
\item[\rm{(2)}] $\bar \kappa_X=0$.
\item[\rm{(3)}] $\bar \alpha_X$ is birational and ${\cal A}_X\setminus\alpha_X(X)$ is finite.
\end{itemize}

\end{theo}

\noindent {\bf  Proof}: we use Theorem~\ref{neut} and the theorem of McQuillan-ElGoul to deduce (2) from (1). To deduce (3) from (2), we use the addition theorem 
of Kawamata for open surfaces and 
Proposition~\ref{neutt}. Finally, we deduce (1) from (3) by elementary methods.\\

\noindent
{\bf (1) implies (2)}:  Suppose that there exists an
algebraically nondegenerate entire curve in $X$ lying in  the foliation defined by a
log $1$-form. By Theorem~\ref{t15}
and Proposition~\ref{nond}
we may assume that $\bar \alpha_X$ is dominant and $\bar \kappa_X\leq 1$. This 
means in particular that there are linearly independent log $1$-forms whose wedge
is not identically $0$. But this wedge gives a log $2$-form and hence 
$\bar \kappa_X\neq -\infty$. But  $\bar \kappa_X\neq 1$ by Theorem~\ref{neut}.
Hence $\bar\kappa_X=0$.\\[-3mm]

\noindent
{\bf (2) implies (3)}:
 Assume that $\alpha_X$ (and 
hence $\bar{\hat \alpha}_X$) is not dominant. Then $\bar{\hat \alpha}_X$
factors through a morphism $\bar \alpha: \bar{\hat X}\r \bar Y\s \bar {\cal A}_X$ 
where $Y=\bar Y\cap {\cal A}_X$ is a hyperbolic curve by Proposition~\ref{nond}
and so $\bar \kappa_Y=1$. 
Let $\bar F$ be a general fiber of $\bar \alpha$, then $\bar F$ is a smooth
projective curve having transversal intersection with $\hat D$, the reduced
total transform of $D$. Hence, by the smoothness of $\bar\alpha$ 
along $\bar F$, we have the following exact sequence
$$0\r \bar K_F\r \Omega(\bar{\hat X},\hat D)|_{\bar F}\r {\cal O}_F\r 0.$$
Hence $\bar K_X|_F=\bigwedge^2\Omega(\bar{\hat X},\hat D)|_{\bar F}=\bar K_F$.
Since $\bar \kappa_X\geq 0$, there is a nontrivial section of some tensor power of 
the log canonical sheaf $\bar K_X$ over $\bar{\hat X}$. Since this section remains
nontrivial over $\bar F$, we see that $\bar \kappa_F\geq 0$. By Theorem~\ref{t10b},
we have $0=\bar\kappa_X\geq \bar \kappa_F+\bar\kappa_B\geq 1$. This
contradiction shows that $\alpha_X$ is dominant (This fact also follows directly
from Theorem 28 of \cite{Kaw}). 
Keeping the same notation as Proposition~\ref{stein} with  the normal surface
$Y=(\bar \gamma)^{-1}({\cal A}_{X})$, then
$\gamma=\bar \gamma|_Y: Y\r  {\cal A}_{X}$ is a surjective finite map and
we have $\bar\kappa_Y\geq 0$ by Theorem~\ref{t10a}.
As $\beta=\bar\beta|_X$ is a morphism to $Y$, we have
$0=\bar\kappa_X\geq \bar\kappa_Y$. Hence, $\bar\kappa_Y=0$
and $Y$ is \'etale over ${\cal A}_{X}$  by Theorem~\ref{t10a}.
This means that $Y$ is a semi-abelian variety and so by the universal
property of the Albanese map $\alpha$, $\gamma:Y\r {\cal A}_{X}$ is an
isomorphism.
By condition (*) and
Proposition~\ref{stein},  we
have that $Y_0=\beta (X)$ is an open subvariety of the normal surface $Y$.
As $\beta=\bar\beta|_X$ is a morphism to $Y_0$, we have
$0=\bar\kappa_X\geq \bar\kappa_{Y_0}\geq \bar\kappa_Y=0$; thus forcing
equality. Hence $Y\setminus Y_0$ is finite by Theorem~\ref{neutt} and the 
result follows. \\[-3mm]

\noindent
{\bf (3) implies (1)}: 
Keeping the same notation as in Proposition~\ref{stein},
$Y_0=\alpha_X(X)$ is an open subset of ${\cal A}_{X}$ with finite
complement.  Let $X^+ =\bar {\hat \alpha}^{-1}(Y_0)$. As $\bar {\hat \alpha}$
is a birational morphism on $X^+$ that restricts to $\alpha_X$ on $X$,
we see that $X^+\setminus X$ is an exceptional divisor of $\bar {\hat \alpha}$
in $X^+$ whose image is contained in the image of the exceptional divisor 
$E$ of $\alpha_X$.
Consider the finite sets $S_1:=\alpha_X(E)\subset {\cal A}_{X} $, 
$S_2:={\cal A}_{X}\setminus Y_0$ and $S'=S_1\cup S_2$. 
 By a translation if necessary, we may assume that the universal covering map 
$h:\cz^2\r  {\cal A}_{X} $ is a morphism 
of additive groups and that $\Gamma=\ker h$ does not intersect $S=h^{-1}(S')$.
Consider the linear map $f_0:\cz \r \cz^2,\ z\mapsto (z, az)$ for $a\in \cz$.
Then $f_0^*w\equiv 0$
for the linear $1$-form $w=dz_2 - a\, dz_1$ where $(z_1,z_2)$ is the standard coordinate
for $\cz^2$.
Now if $f_1=h\circ f_0$ is algebraically
degenerate, then its
image lies in an elliptic curve or a rational curve in $\bar {\cal A}_{X} $
and the latter must intersect the boundary of ${\cal A}_{X} $ at at least 
and hence exactly two points. In either case $f_1$ is an \'etale 
covering over its image, which is either $\cz^*$ or an elliptic curve, 
and so $f_0^{-1}(\Gamma)=\ker f_1\neq \{(0,0)\}$. Therefore, 
since log $1$-forms on $ {\cal A}_{X} $ correspond bijectively with linear
$1$-forms on $\cz^2$, it suffices to produce an $a$ such that  
$f_0^{-1}(\Gamma'\cup S)=\emptyset$, $\Gamma'=\Gamma\setminus \{(0,0)\}$,
for then $f_1(\cz)$ does not
intersect $S'$ and $\alpha_X^{-1}\circ f_1$ would be the required holomorphic curve
on $X$ as $\alpha_X^{-1}$ is a well defined holomorphic map outside $S'$.
For this, we only need to choose $a\in \cz$ outside the countable set
$\{v/u\ |\ (u,v)\in \Gamma'\cup S\}$.
\qed

 \vspace{1cm}
 The 
additional condition (*) on $\alpha_X$ 
is essential for Theorems~\ref{neut} and \ref{neupt},
as the following counterexamples shows. We remark here however
that it can be weakened still, especially for Theorem~\ref{neupt}, as
can be seen from the proofs we give. We also remark that one can prove that (2) and (3) of Theorem~\ref{neupt} are equivalent even without condition (*).

\begin{prop}\label{p6}
(Counterexample for $\bar \kappa_{X} =1$)
Let $E$ be an elliptic curve and $p:E \times E \rightarrow E$
the projection to the first factor. Let $P_{1}, P_{2} \in E$ be two distinct 
points, and $Q_{i} \in p^{{-1}}(P_{i})$, $i=1,2$ two points. 
Let $b: \bar X \rightarrow E \times E$
be the blow up of $E \times E$ in the points $Q_{1}, Q_{2}$. Let $D$ be the union
of the proper transforms of $p^{{-1}}(P_{i})$, $i=1,2$ in $\bar X$, and 
$X := \bar X \setminus D$. Then $\bar q_X = 2$ and $\bar \kappa_{X}=1$,
but $X$ admits entire curves $f:\cz \rightarrow X$ which are neither algebraically
degenerate nor linearly degenerate (i.e. the condition given in \rm{(1)} of
Theorem~\ref{neupt}). 
\end{prop}
\noindent {\bf Proof}:
First it is easy to see that $\bar q_{X} \geq 2$, since 
$q_{E \times E}=2$ and linearly independent 1-forms on $E \times E$ lift to 
linearly independent 1-forms on $\bar X$.

The fact that there exist entire curves $f:\cz \rightarrow X$ which are not
algebraically degenerate is an easy application of the main result of 
Buzzard-Lu'00 \cite{BL}: The map $p \circ b : X \rightarrow E$ is a surjective
holomorphic map defining $X$ into an elliptic fibered surface over the curve $E$ 
without any branching. So by Theorem 1.7 of \cite{BL}, $X$ is dominable by $\cz^{2}$
and there are (a lot of) entire curves $f:\cz \rightarrow X$ which are neither
algebraically degenerate nor linearly degenerate 
(i.e. the condition given in \rm{(1)} of
Theorem~\ref{neupt}). By Theorem~\ref{t1} this still implies in particular
that $\bar q_{X} \leq 2$, so $\bar q_{X} = 2$.

It is easy to see that $\bar \kappa_{X}\leq 1$: Since $K_{E \times E}=0$
and by the standard formula for the relation of the canonical divisors under blow
ups of points we get 
$$\bar K_{X} = (p \circ b)^{-1}(P_{1}) + (p \circ b)^{-1}(P_{2})\, .$$ 
Hence, for a generic point $P \in E$, the fiber $F:= (p \circ b)^{-1}(P)$
has the property that the restriction to $F$ of the line bundles 
$(\bar K_{X})^{\otimes m}$ are trivial, hence, their sections cannot seperate points.
Hence, $\bar \kappa_{X} \leq1$.

It remains to prove that $\bar \kappa_{X} \geq 1$. By the Riemann-Roch theorem
for the curve $E$ and the divisor $P_{1}+P_{2} \subset E$ of degree $2$ (see
for example Hartshorne '77 \cite{H}, page 295), we get that for the log curve
$E_{0}:= E \setminus \{P_{1}, P_{2}\} \subset E$, we have
$$h^{0}(E,  \bar K_{E_{0}})= h^{0}(E,  P_{1}+P_{2})= {\rm deg}(P_{1}+P_{2})= 2\, .$$
Hence, there exist two linearly independent log 1-forms $\eta_{1}, \eta_{2}$
on $E_{0} \subset E$. The key point is now that we can use them to construct
two linearly independent sections of the line bundle $\bar K_{X}$ on $\bar X$
(since this will prove $\bar \kappa_{X} \geq 1$):

 Locally, let $x$ be the base coordinate and $y$ the fiber 
coordinate of the projection $p:E \times E \rightarrow E$, 
both linear with respect to the linear structure of $E$.  We blow up the origin
of $(x,y)$ and let $X= \bar X \setminus D$, where $D$ is the proper transform of $x=0$ by $b$. 
Then
in a neighborhood of the point of interection of the exceptional curve of $b$ and the
proper transform of $y=0$ by $b$, $X$ is parametrized by $x,z$ where $y=zx$. Any log 1-form
of $E\setminus \{0\}$ is of the form $f(x)dx/x$ near $0$. Now 
$$dy=zdx+xdz$$ is a 1-form on $\bar X$ and, hence, 
\begin{equation}
f(x){dx\over x}\wedge dy= f(x) dx\wedge dz \label{form}
\end{equation}
is a log 2-form on 
$(\bar X,  \{x=0\})$, which by (\ref{form}) doesn't have any poles  on $X$
and, hence, is  a log 2-form on $(\bar X, D)$. By this local argument we thus see 
that the forms $\eta_{1}\wedge dy, \eta_{2}\wedge dy$ are
two linearly independent global sections of the line bundle $\bar K_{X}$ on $\bar X$. \qed

   \noindent Gerd Dethloff\\
   Universit\'{e} de Bretagne Occidentale\\
    UFR Sciences et Techniques \\
D\'{e}partement de Math\'{e}matiques\\
6, avenue Le Gorgeu, BP 452 \\
    29275 Brest Cedex, France \\
    e-mail: gerd.dethloff@univ-brest.fr\\

    \noindent  Steven S.-Y. Lu\\
    D\'epartement de Math.\\
    Universit\'{e} du Quebec \`a Montreal\\
    201 Ave. du Pr\'esident Kennedy\\
    Montr\'eal H2X 3Y7\\
    Canada\\
    e-mail: lu@math.uqam.ca

\end{document}